\documentclass[11pt]{amsart}

\usepackage{amssymb,amsmath,amsthm}
\usepackage{stmaryrd}

\usepackage{diagrams}
\diagramstyle{PostScript=dvips,nohug,nopixel,amstex,labelstyle=\scriptstyle}
\newarrow{To}----{->}
\newarrow{Onto}----{->>}
\newarrow{Into}{C}---{->}

\newcommand{\Rp}{R_{+}}
\newcommand{\Rm}{R_{-}} 
\newcommand{\Rpm}{R_{\pm}}

\newcommand{\Rsc}{R_{\mathrm{sc}}}
\newcommand{\Rac}{R_{\mathrm{ac}}}

\newcommand{\SRR}{SR_{\R}}

\newcommand{\tensor}{\mathbin{\widetilde{\otimes}}}

\newcommand{\dirac}{\mbox{$\not\negthinspace\partial$}}
\newcommand{\Z}{\mathbb{Z}}
\newcommand{\R}{\mathbb{R}}
\newcommand{\C}{\mathbb{C}}
\newcommand{\g}{\mathfrak{g}}
\DeclareMathOperator{\ad}{ad}
\DeclareMathOperator{\Cl}{\mathbb{C}l}

\DeclareMathOperator{\End}{End}

\theoremstyle{plain}
\newtheorem{theorem}{Theorem}[section]
\newtheorem{proposition}[theorem]{Proposition}
\newtheorem{lemma}[theorem]{Lemma}
\newtheorem*{schur}{Schur's Lemma}
\newtheorem{corollary}[theorem]{Corollary}
\theoremstyle{remark}
\newtheorem{remark}[theorem]{Remark}
\newtheorem*{example}{Example}
\theoremstyle{definition}
\newtheorem*{definition}{Definition}

\title{Representation rings of Lie superalgebras}
\author{Gregory D. Landweber}

\email{greg@math.uoregon.edu}
\urladdr{http://www.uoregon.edu/\~{}greg/}
\address{Mathematics Department\\
         University of Oregon\\
         Eugene, OR 97403-1222}

\begin{document}


\begin{abstract}
Given a Lie superalgebra $\g$, we introduce several variants of the representation ring, built as subrings and quotients of the ring $R_{\Z_{2}}(\g)$ of virtual $\g$-supermodules, up to (even) isomorphisms. In particular, we consider the ideal $R_{+}(\g)$ of virtual $\g$-supermodules isomorphic to their own parity reversals, as well as an equivariant $K$-theoretic super representation ring $SR(\g)$ on which the parity reversal operator takes the class of a virtual $\g$-supermodule to its negative.  We also construct representation groups built from ungraded $\g$-modules, as well as degree-shifted representation groups using Clifford modules. The full super representation ring $SR^{\bullet}(\g)$, including all degree shifts, is then a $\Z_{2}$-graded ring in the complex case and a $\Z_{8}$-graded ring in the real case.

Our primary result is a six-term periodic exact sequence relating the rings $R^{\bullet}_{\Z_{2}}(\g)$, $R^{\bullet}_{+}(\g)$, and $SR^{\bullet}(\g)$. We first establish a version of it working over an arbitrary (not necessarily algebraically closed) field of characteristic $0$. In the complex case, this six-term periodic long exact sequence splits into two three-term sequences, which gives us additional insight into the structure of the complex super representation ring $SR^{\bullet}(\g)$. In the real case, we obtain the expected 24-term version, as well as a surprising six-term version, of this periodic exact sequence.
\end{abstract}

\maketitle


%

\tableofcontents

\section{Introduction}

This paper investigates several possible definitions for the representation ring of a Lie superalgebra and delves into the structure of these various rings (and ideals and groups) by examining maps between them. This work originated from the author's attempt to reconcile a version of the Grothendieck group of virtual representations used by some superalgebraists \cite{BK0} with an equivariant $K$-theoretic approach used by topologists \cite{FHT}. A Lie superalgebra is the $\Z_{2}$-graded generalization of a Lie algebra, with additional signs consistently interspersed in the definition. When working with representations of Lie superalgebras, one usually considers $\Z_{2}$-graded representations on super vector spaces, which decompose into two pieces
$V = V_{0} \oplus V_{1}$ according to the grading. Given such a super vector space $V$, we can reverse its parity by interchanging its even and odd degree components, but leaving all other structures unchanged. We call this parity reversed supermodule $V^{\Pi} = V_{1} \oplus V_{0}$. In some papers, such as \cite{BK0}, the class of a supermodule $V$ is identified with that of its parity reversal $V^{\Pi}$ in the representation ring. In other words, this representation ring is constructed so that the parity reversal operator $\Pi$ acts as the identity.

On the other hand, recall that the Grothendieck $K$-construction transforms the abelian semi-ring of isomorphism classes of vector bundles into a ring by introducing formal differences. This $K$-theory is then the natural setting in which to consider the index
$$\mathrm{Index}\,D = \mathrm{Ker}\,D - \mathrm{Im}\,D$$
of (families of) Fredholm operators $D$. Going in the opposite direction, given a super vector space $V = V_{0} \oplus V_{1}$, constructions involving formal differences $V_{0} - V_{1}$ can often be expressed in terms of the index of an odd linear operator $D : V \to V$ which interchanges the two components. An example of this is the Weyl character formula,
which expresses the character of the irreducible $\g$-module $V_{\lambda}$ with highest weight $\lambda$ as
\begin{equation}\label{eq:weyl}
( V_{\lambda} \otimes \mathbb{S}_{0} ) -  ( V_{\lambda} \otimes \mathbb{S}_{1} )
   = \sum_{w\in W_{\g}} (-1)^{w}\, e^{w(\lambda+\rho)},
\end{equation}
where the computation is performed in the representation ring $R(\mathfrak{t})$
of the Cartan subalgebra $\mathfrak{t}$. The left side involves a formal difference of the two components of the $\Z_{2}$-graded spin representation $\mathbb{S} = \mathbb{S}_{0} \oplus \mathbb{S}_{1}$ of the Clifford algebra $\mathrm{Cl}(\g/\mathfrak{t})$, viewed as a superalgebra. In his beautiful paper \cite{B}, Bott suggests viewing the Weyl character formula (\ref{eq:weyl}) (as well as the related Borel-Weil-Bott theorem) in terms of the index of a Dirac operator on the coset space
$G/T$. This view is espoused in \cite{Ko} by Kostant, who expresses a generalization \cite{GKRS} of the Weyl character formula (\ref{eq:weyl}), replacing the Cartan subalgebra $\mathfrak{t}$ with a general subalgebra of maximal rank, in terms of the index of a formal algebraic Dirac operator $\dirac : V_{\lambda} \otimes \mathbb{S}_{0} \to V_{\lambda} \otimes \mathbb{S}_{1}$. In addition, the author's generalization of this result \cite{La} to the Kac-Moody setting has been recently used in a construction of twisted equivariant $K$-theory \cite{FHT} by Freed, Hopkins, and Teleman.

One of our goals here is to introduce and study an equivariant $K$-theoretic construction, which we will call the super representation ring (in analogy to the supertrace), that is essentially built from formal differences $V_{0} - V_{1}$ of the two components of supermodules $V = V_{0} \oplus V_{1}$. In particular, we want the parity reversal operator, when acting on the super representation ring, to take the class of a representation to its negative: $[V^{\Pi}] = -[V]$. (We note that our construction of the super representation ring has a similar spirit to that of the variant $K_{\pm}(X)$ for spaces $X$ with involution recently introduced by Atiyah and Hopkins in \cite{AH}.) If we compare how these two rings behave for the complex Clifford algebra $\Cl(n)$, we will find that
the construction identifing supermodules with their parity reversals
gives $\Z$ for all $n$, while our super representation ring alternates between $\Z$ and $0$ in even and odd degrees, respectively, reproducing the $K$-theory of a point.  Likewise, for real Clifford algbras, our super representation ring gives the $KO$-theory of a point, and the real case possesses a much richer structure. Extending this $K$-theoretic connection, we construct degree-shifted versions of our super representation ring by considering supermodules carrying additional Clifford actions.
(In \cite{FHT} and \cite{FHT1}, Freed, Hopkins, and Teleman refer to the odd actions of the Clifford generators as ``supersymmetries'', in the sense that they are automorphisms of the supermodules which interchange the even and odd components.)
Considering all degrees simultaneously, in the complex case this gives us a $\Z_{2}$-graded ring which contains information similar to that of the
representation ring constructed by identifying supermodules with their parity reversals, while in the real case we obtain a $\Z_{8}$-graded ring closely related to equivariant $KO$-theory.

Our main results are Theorems \ref{exact-sequence1} and \ref{24}, which give periodic long exact sequences relating the homogeneous components of these representation rings:
\begin{equation*}
	\cdots \to 
	R_{\Z_{2}}^{-n-1}(\g) \xrightarrow{i^{*}} \Rp^{-n}(\g) \to
	SR^{-n}(\g)
	\to R_{\Z_{2}}^{-n}(\g) \xrightarrow{i^{*}} \Rp^{-n+1}(\g)\to \cdots
\end{equation*}
where $SR^{\bullet}(\g)$ is our graded super representation ring, $R_{\Z_{2}}^{\bullet}(\g)$ is the graded ring of virtual supermodules, and $R_{+}^{\bullet}(\g)$ is the graded ring (or more precisely, ideal) of virtual supermodules satisfying $V \cong V^{\Pi}$.
The degree zero component of the self-dual representation ring $R_{+}^{0}(\g)$ is isomorphic to the alternative version of the representation ring constructed by identifying supermodules with their parity reversals, although this is an isomorphism only of additive groups, not of rings.
We first establish a six-term version of this periodic long exact sequence in a slightly different form in Theorem \ref{exact-sequence}, where we work over an arbitrary field of characteristic 0 which is not necessarily algebraically closed. There we perform the bulk of our computations, and we employ both graded and ungraded representations in place of degree shifts. We then specialize to the complex case, where the corresponding six-term periodic exact sequence of Theorem \ref{exact-sequence1} splits into two short exact sequences, and much of the discussion simplifies significantly.  This can be viewed as a warmup for the real case, in which we obtain the expected 24-term periodic long exact sequence given in Theorem \ref{24}.  However, we were surprised to find that the six-term periodic long exact sequence of Theorem \ref{exact-sequence} also generalizes in a different way to give a collection of eight six-term periodic long exact sequences in the real case!

Our use of Clifford algebras to build the homogeneous components of the graded representation rings is motivated by the approach to $K$-theory introduced by Atiyah, Bott, and Shapiro in \cite{ABS} and examined extensively
by Karoubi in \cite{Ka}.  In Section \ref{subsection-representation-ring}, we review the Atiyah-Bott-Shapiro construction of the $K$-theory of a point and generalize it to the equivariant case of the degree-shifted representation modules for a compact Lie group. (We work with a Lie group rather than a Lie algebra as a nod to the geometric nature of \cite{ABS}.) In Section \ref{lie-superalgebras} we shift to the algebraic setting, introducing Lie superalgebras and defining our various representation rings and groups.  In Section \ref{involutions} we investigate the parity reversal involution as well as a related conjugation involution for ungraded modules, and we use these to construct the first version of our periodic exact sequence in Section \ref{les} by examining the kernels and cokernels of the maps used to define the super representation ring.  These two sections go into significant detail for the following two reasons: First, we are working over an arbitrary field of characteristic 0, which allows us to explore the real case, which has a much richer structure than the complex case. Second, this treatment is intended for an audience interested in the topological aspects and inexperienced with Lie superalgebras.  Experts are urged to skim or even skip these sections and go directly to Section \ref{degree-shift}, where we reintroduce Clifford algebras and construct the degree-shifted representation groups, which we use to prove the complex version of our periodic long exact sequence in Section \ref{redux}. Finally, in Section \ref{real} we consider real representation rings and briefly outline how the results of all the previous sections generalize, giving both the 24-term and six-term versions of our periodic exact sequence.

In \cite{La2}, we go on to discuss representation rings of Lie supergroups and connections with Dirac operators 
and twisted $K$-theory. There, we consider the restriction and induction maps on twisted versions of the super representation ring.
For a taste of this, let $G$ be a compact, simply connected Lie group, and consider the Lie supergroup $\Pi(T^{*}G)$, the parity reversal of the cotangent bundle $T^{*}G$, which has $G$ as its underlying even part and the parity reversal of the coadjoint representation $\g^{*}$ for its odd directions.  The corresponding Lie superalgebra is $\g \oplus \Pi\g^{*}$.  Let $T$ be a maximal
torus inside $G$. Consider the following commutative diagram, which relates the twisted super representation rings of the Lie supergroups with the standard representation rings of the Lie groups:
$$\begin{diagram}[tight,width=3.5em]
SR^{\tau}( \Pi T^{*}G ) & \rTo^{\cong} & R(G) \\
\dTo & & \dTo \\ 
SR^{\tau}( \Pi T^{*}T ) & \rTo^{\cong} & R(T)
\end{diagram}$$
where the top and bottom maps are Thom isomorphisms and the map on the left side is the restriction map.
The right side then gives the Weyl character formula, mapping an irreducible $G$-module $V_{\lambda}$ to the expression (\ref{eq:weyl}) in $R(T)$.  Going in the opposite direction, we can also consider the induction map on representation rings, as described by Bott in
\cite{B0}, and we obtain the commutative diagram
$$\begin{diagram}[tight,width=3.5em]
SR^{\tau}( \Pi T^{*}G ) & \lTo^{\cong} & R(G) \\
\uTo & & \uTo \\ 
SR^{\tau}( \Pi T^{*}T ) & \lTo^{\cong} & R(T)
\end{diagram}$$
where once again the top and bottom map are Thom isomorphisms, and the left side is the standard induction map. The right side is then Dirac induction, taking a one-dimensional $T$-module $\C_{\mu}$ with weight $\mu$ to the $G$-equivariant index of the Dirac operator on the homogeneous space $G/T$ with values in the homogeneous bundle induced by $\C_{\mu}$:
$$\C_{\mu} \mapsto \mathrm{Index}\,\dirac_{G/T} \otimes \C_{\mu} = (-1)^{w}\,V_{w\mu - \rho}.$$
If we replace $T$ with a subgroup of $G$ containing $T$,
then this method reproduces the generalization \cite{GKRS} of the Weyl character formula, as well as the Dirac induction map of \cite{La0}. In the infinite dimensional setting, such Dirac induction maps are a fundamental ingredient in understanding
the Verlinde algebra via twisted equivariant $K$-theory, as described in \cite{Fr}.  In addition, the author expects constructions involving the real version of the super representation ring to play a similar role in twisted, equivariant $KO$-theory.

Dimitry Leites and the referee have both pointed out that the representation rings considered here in fact describe the ``body'' of a superring. Our discussion deals only with the geometric points of representations, but in the presence of odd parameters, one can consider additional representations, such as ones involving the queerdeterminant (see \cite[\S 1.8]{SoS}). The full representation superring will not be considered here, but it is an interesting project for future study.

\smallskip
\noindent
\textbf{Acknowledgements.}
The author would like to thank his colleagues Jon Brundan and Dev Sinha for their insight and helpful conversations. The commutative diagrams in this article were typeset using the ``diagrams''  \TeX\ macro package by Paul Taylor.

\section{The representation ring of a compact Lie group}\label{subsection-representation-ring}

Let $G$ be a compact Lie group. Recall that its representation ring $R(G)$ is the Grothendieck ring of isomorphism classes of finite dimensional, complex $G$-modules. (By complete reducibility, $R(G)$ is furthermore the free $\Z$-module generated by the isomorphism classes of irreducible, finite dimensional, complex representations of $G$.) Viewing such representations as $G$-equivariant complex vector spaces, we can alternatively define the representation ring $R(G)$ to be $K_{G}(\mathrm{pt})$, the $G$-equivariant $K$-theory of a point. Taking this point of view, we can use an equivariant version of the Atiyah-Bott-Shapiro isomorphism \cite{ABS} to define a degree-shifted representation module
$R^{-n}(G) = K_{G}^{-n}(\mathrm{pt}) = \tilde{K}_{G}(S^{n})$ as follows: 

Let $\Cl(n)$ denote the complex Clifford algebra on the $n$ generators
$e_{1},\ldots,e_{n}$ with the relations
\begin{equation}\label{eq:clifford}
e_{i}^{2} = -1 \qquad \text{and} \qquad e_{i} \cdot e_{j} = - e_{j} \cdot e_{i} \text{ for $i\neq j$}.
\end{equation}
The Clifford algebra is a $\Z_{2}$-graded algebra, or superalgebra,
decomposing into its even and odd degree parts:
$\Cl(n) = \Cl_{0}(n) \oplus \Cl_{1}(n)$. In the supersymmetric world,
everything carries a $\Z_{2}$-grading, so we consider representations
in which the Clifford algebra acts on a $\Z_{2}$-graded vector space,
or super vector space, $V = V_{0}\oplus V_{1}$, respecting the gradings.
In other words, the action of the even part
$\Cl_{0}(n)$ preserves the grading on $V$, while the action of
the odd part $\Cl_{1}(n)$ interchanges $V_{0}$ and $V_{1}$. Let $M(n)$
denote the Grothendieck group of isomorphism classes of such finite dimensional $\Z_{2}$-graded
$\Cl(n)$-modules. (Here, we require that isomorphisms of $\Z_{2}$-graded modules be even maps which preserve the grading unless otherwise noted.) The inclusion $i(n) : \Cl(n) \hookrightarrow \Cl(n+1)$
of Clifford algebras induces a restriction map
$$i(n)^{*} : M(n+1) \rightarrow M(n)$$ for $\Z_{2}$-graded Clifford modules.
Atiyah, Bott, and Shapiro proved in \cite{ABS} that
\begin{equation}\label{eq:ABS}
K^{-n}(\mathrm{pt}) \cong \mathrm{Coker}\, i(n)^{*} = M(n) \, / \, \mathrm{Im}\,i(n)^{*}
\end{equation}
as abelian groups.

Repeating this construction in the equivariant context, let $M_{G}(n)$
denote the Grothendieck group of isomorphism
classes of finite dimensional $\Z_{2}$-graded $G\times\Cl(n)$-modules, i.e.,
Clifford modules possessing an even ($\Z_{2}$-degree preserving) $G$-action commuting with the $\Cl(n)$-action (more precisely, we are considering representations of the tensor product $\C[G] \otimes \Cl(n)$, where $\C[G]$ is the complex group ring of $G$). Once again we have a restriction map
$$i_{G}(n)^{*} : M_{G}(n+1) \rightarrow M_{G}(n)$$ for these equivariant $\Z_{2}$-graded Clifford modules, and following the construction of Atiyah-Bott-Shapiro, we make the following:

\begin{definition}\label{Rminusn}
The degree-shifted representation module $R^{-n}(G)$ is defined by
\begin{equation}\label{eq:Rminusn}
R^{-n}(G) := \mathrm{Coker}\, i_{G}(n)^{*} = M_{G}(n) \, / \, 	\mathrm{Im}\,i_{G}(n)^{*},
\end{equation}
the classes of virtual, finite dimensional, complex $G$-super\-mod\-ules with $n$ supersymmetries, modulo those admitting an additional supersymmetry.
\end{definition}

Each of these degree-shifted groups $R^{-n}(G)$ is in fact a module over the standard representation ring $R^{0}(G) \cong R(G)$, with action given by the tensor product. To clarify the second half of this definition, which is in the language of \cite{FHT} and the appendix of \cite{FHT1}, note that the Clifford generators $e_{j}$ act on a $\Z_{2}$-graded Clifford module $V$ by odd automorphisms, interchanging the $V_{0}$ and $V_{1}$ components. Such odd automorphisms are often referred to as supersymmetries.  We say that two homogeneous $\Z_{2}$-graded objects $x$ and $y$ supercommute if $x\cdot y = (-1)^{|x|\,|y|}\,y\cdot x$, where $|x|$ and $|y|$ denote their even or odd $\Z_{2}$-degrees.
In this case, the generators $e_{j}$ have odd degree and anti-commute due to the Clifford relations (\ref{eq:clifford}). Note that each of these supersymmetries does not supercommute with itself, which for an odd object $x$ would require $x\cdot x = -x\cdot x$ and therefore $x^{2} = 0$. By the Clifford relations (\ref{eq:clifford}), each of these generators actually squares to $-1\,\mathrm{Id}$, but since we are working over $\C$, we could just as well replace them with the generators $i\,e_{j}$ which square to $\mathrm{Id}$.

Note that these degree-shifted representation modules do not contain any additional information about the Lie group. Since we are considering $G\times\Cl(n)$-modules, where the $G$-action commutes with the Clifford action, any such module can be built as the direct sum of
tensor products $V \otimes \mathbb{S}$, where $V$ is a $G$-module
and $\mathbb{S}$ is a $\Z_{2}$-graded Clifford module. We therefore
have
\begin{equation}\label{eq:rn}
	M_{G}(n) \cong R(G)\otimes M(n),
\end{equation}
and the restriction map becomes $i_{G}(n)^{*} \cong 1 \otimes i(n)^{*}$. In light of the isomorphism (\ref{eq:ABS}), our definition (\ref{eq:Rminusn}) gives us for non-negative integers $n$ an isomorphism of $R(G)$-modules:
\begin{equation*}
R^{-n}(G) \cong R(G) \otimes K^{-n}(\mathrm{pt})
\cong \begin{cases}
	R(G) & \text{ if $n$ is even,} \\
	0    & \text{ if $n$ is odd.}
\end{cases}
\end{equation*}
The even/odd decomposition of $K^{n}(\mathrm{pt})$ comes courtesy of Bott periodicity (or more precisely the 2-fold periodicity of complex Clifford algebras
that underlies Bott periodicity).
By Bott periodicity, we have an isomorphism $R^{-n}(G) \cong R^{-n-2}(G)$,
which allows us to define $R^{n}(G)$ for all integers $n$.
Incorporating these degree shifts, $R^{\bullet}(G)$ is then a $\Z_{2}$-graded ring extending the usual representation ring $R(G) = R^{0}(G)$, and it is isomorphic to $R(G) \otimes K^{\bullet}(\mathrm{pt})$ as a graded ring over $R(G)$. The product structure
on $R^{\bullet}(G)$ is given by the tensor product of Clifford modules; if $V$ and $W$ are $\Cl(n)$ and $\Cl(m)$-supermodules, respectively, then their tensor product $V \otimes W$ is a $\Cl(n+m)$-supermodule.

Since we are working with $\Z_{2}$-graded modules, we can introduce
the parity reversal operator $\Pi$. Given a supermodule $V$,
we define the parity reversed supermodule $V^{\Pi}$ to be the same underlying vector space with the same action as $V$, but with the opposite grading: $(V^{\Pi})_{0} = V_{1}$ and $(V^{\Pi})_{1} = V_{0}$. (In particular, $V^{\Pi}$ is the right parity reversal, and unlike the more conventional left parity reversal $\Pi V$, we do not introduce any signs into the action. See Remark~\ref{remark:pi} below.) Parity reversal gives an involution on the isomorphism classes of virtual supermodules, $\Pi : M_{G}(n) \to M_{G}(n)$.
 
In the physics literature, the grading on
a supermodule $V$ is specified by the operator $(-1)^{F}$, which acts by $+1$ on $V_{0}$ and by $-1$ on $V_{1}$. Here, the ``fermion number'' $F$ is simply the even or odd $\Z_{2}$-degree. (For $\Z$-graded supermodules, such as exterior algebras, the fermion number $F$ can actually be $\Z$-valued, but the operator $(-1)^{F}$ still gives only a $\Z_{2}$-grading.) Even elements are called ``bosons'', while odd elements are called ``fermions''. The parity reversal operator therefore interchanges bosons with fermions, and this duality is known to physicists \cite{G} as a ``Klein flip''.
This grading operator $(-1)^{F}$ commutes with even homomorphisms, such as the $G$-action, and anti-commutes with odd homomorphisms, such as the actions of the Clifford algebra generators.  If $V$ has grading
operator $(-1)^{F}$, then $V^{\Pi}$ has grading operator $-(-1)^{F}$.

Our definition (\ref{eq:Rminusn}) can be stated quite elegantly in terms of the parity reversal operator.  In (\ref{eq:Rminusn}), the subgroup by which we quotient consists of supermodules isomorphic to their own parity reversal. However, since parity reversal is an involution on the space of supermodules, we attempt to quantize this relation by restricting the possible isomorphisms to those which are themselves involutions.

\begin{proposition}\label{prop-1}
The image $\mathrm{Im}\,i_{G}(n)^{*}$ in (\ref{eq:Rminusn}) is generated by classes of $\Z_{2}$-graded Clifford modules $[V]\in M_{G}(n)$ which admit 
a $G\times \Cl(n)$-equivariant parity reversing involution $\alpha : V \to V^{\Pi}$.
\end{proposition}

(More precisely, a $G\times\Cl(n)$-equivariant isomorphism $\alpha: V\to V^{\Pi}$ induces a parity reversed isomorphism $\alpha^{\Pi} : V^{\Pi} \to (V^{\Pi})^{\Pi} = V$ with the same action as $\alpha$ on the underlying Clifford module. So, the condition
in the Proposition can be stated as $\alpha^{\Pi} \circ \alpha = \mathrm{Id} : V\to V$, or by a slight abuse of notation, simply $\alpha^{2} = \mathrm{Id}$.)

\begin{proof}
In light of (\ref{eq:rn}), this is really a statement about the image of $i(n)^{*}$, as the $G$-action commutes with everything.
If $[V] \in \mathrm{Im}\,i_{G}(n)^{*}$, then $V$ admits an action of $e_{n+1}$ satisfying the Clifford relations (\ref{eq:clifford}).
Let $\alpha = e_{n+1}\,(-1)^{F}$. Then $\alpha : V\to V$ is an
odd automorphism, which can be viewed as an (even) isomorphism
$\alpha : V \to V^{\Pi}$. Furthermore, due to the added $(-1)^{F}$ factor, this $\alpha$ now commutes (in the usual
non-supersymmetric sense) with the actions of each $e_{j}$ for $j\leq n$, and we have $\alpha^{2} = \mathrm{Id}$. Conversely, given a class $[V]\in M_{G}(n)$ with an equivariant isomorphism $\alpha : V \to V^{\Pi}$ satisfying $\alpha^{2} = \mathrm{Id}$, we can set $e_{n+1} = \alpha\,(-1)^{F}$. The action of this $e_{n+1}$
then satisfies the Clifford relations (\ref{eq:clifford}).
\end{proof}

It follows that the parity reversal operator $\Pi : M_{G}(n) \to M_{G}(n)$ acts trivially on all elements in the image $\mathrm{Im}\,i_{G}(n)^{*}$, and thus it
descends to an involution of the degree-shifted representation module $R^{-n}(G)$. Given any $G\times \Cl(n)$-supermodule $V$, the supermodule $V \oplus V^{\Pi}$ clearly admits a parity reversing involution, and thus $[V^{\Pi}] = - [V]$ in $R^{-n}(G)$.  When working over the complex numbers, we have the following supersymmetric version of Schur's Lemma (see \cite{Kac1}, \cite{Kac}, or
\cite[\S 2]{BK1}):

\begin{schur}
Let $A$ be a collection of operators (both even and odd) acting irreducibly on a super vector space $V = V_{0} \oplus V_{1}$ over an algebraically closed field. The only even $A$-equivariant endomorphisms of $V$ are scalar multiples of the identity, and there are two possibilities for the odd endomorphisms: either there are no nonzero odd endomorphisms, or the odd endomorphisms consist of scalar multiples of an $A$-equivariant parity reversing involution of $V$.
\end{schur}

In light of Schur's Lemma, every complex Clifford module which is isomorphic to its own parity reversal must admit a parity reversing involution. However, for real representations and $KO$-theory, 
there is a richer structure.  In the real case, there are indeed
Clifford modules which are isomorphic to their own parity reversals but
which do not admit a parity reversing involution. For example, when constructing
$KO^{-1}(\mathrm{pt})$, we start with $M(1) \cong \Z$, which is generated
under direct sums by the class of the super vector space $V = \R \oplus \R$
with the Clifford action
$$e \mapsto \begin{pmatrix} 0 & -1 \\ 1 & \phantom{-}0 \end{pmatrix}.$$
By a simple dimension count, we must then have $V^{\Pi} \cong V$, and thus
every $\mathrm{Cl}(1)$-module is isomorphic to its own parity reversal.
Here, the $\mathrm{Cl}(1)$-equivariant isomorphism is none other than the action of $e$ itself. However, since $e^{2} = -1$, this isomorphism is not an involution. In
the complex case, we could simply multiply the action of $e$ by $\sqrt{-1}$, but we do not have
that luxury here. In order to obtain an honest involution, we must instead work with
$2V = V \oplus V \cong V \oplus V^{\Pi} = \R^{2}\oplus \R^{2}$, which gives us enough room to construct
a square root of $-\mathrm{Id}$. So, we have $R\R^{-1}(\R) \cong KO^{-1}(\mathrm{pt}) \cong \Z_{2}$. In general, the degree-shifted real representation modules are given by
\begin{equation*}
R\R^{-n}(G) \cong R\R(G) \otimes KO^{-n}(\mathrm{pt})
\cong \begin{cases}
	R\R(G) & \text{ if $n \equiv0,4\bmod 8$,} \\
	R\R(G) \otimes \Z_{2} & \text{ if $n\equiv 1,2\bmod 8$,} \\
	0 & \text{ otherwise,}
\end{cases}
\end{equation*}
and indeed we have an isomorphism $R\R^{\bullet}(G) \cong R\R(G) \otimes KO^{\bullet}(\mathrm{pt})$ of graded rings over the real representation ring $R\R(G)$. We will return to the real case in Section \ref{real}.

The above discussion applies equally well if we consider representations of a compact Lie algebra $\g$ in place of the compact Lie group $G$. In the Lie algebra case, we can define the full representation ring $R^{\bullet}(\g)$,
which is a $\Z_{2}$-graded ring over the Grothendieck ring $R(\g) \cong R^{0}(\g)$ of isomorphism classes of finite dimensional, complex $\g$-modules. If $G$ is a compact, simply connected Lie group with Lie algebra $\g$, then $R^{\bullet}(G) \cong R^{\bullet}(\g)$.

\section{Lie superalgebras}\label{lie-superalgebras}

We next consider representation rings for Lie superalgebras. Here we work with Lie superalgebras and representations over a field $k$ of characteristic $0$ (we do not assume that $k$ is algebraically closed), with an eye towards later sections in which we specify $k=\C$ or $k=\R$. A Lie superalgebra (see, for example, \cite{Kac1}, \cite{Kac}) is a
$\Z_{2}$-graded vector space $\g = \g_{0}\oplus\g_{1}$
along with a $\Z_{2}$-graded bracket
$[\cdot,\cdot]\colon \g\otimes\g\to\g$ which is anti-symmetric
in the supersymmetric sense. In other words, for homogeneous elements $X,Y\in \g$ with $\Z_{2}$-degrees $|X|,|Y|$,
respectively,
we have
$$[X,Y]\in \g_{|X|+|Y|}, \quad \text{and} \quad
[X,Y] = - (-1)^{|X|\,|Y|}[Y,X].$$
The Jacobi identity for
Lie superalgebras is best stated as the fact that the adjoint map
$\ad_{X} = [X,\cdot\,]$ is a super-derivation, i.e.,
\begin{equation}\label{eq:jacobi}
	\bigl[X, [Y,Z] \bigr]
	= \bigl[[X,Y],Z\bigr]
	+ (-1)^{|X|\,|Y|}\bigl[Y,[X,Z]\bigr].
\end{equation}
(In general, a Lie superalgebra behaves like a Lie
algebra except that interchanging two neighboring odd
elements in any identity incurs a factor of $-1$.)
Alternatively, a
Lie superalgebra consists of a standard Lie algebra
$\g_{0}$ for the even part and a representation
of $\g_{0}$ on a module $\g_{1}$ for the odd part, equipped with a 
$\g_{0}$-invariant symmetric bilinear form $\g_{1}\otimes\g_{1}\to\g_{0}$, furthermore satisfying the Jacobi identity
$$ \bigl[ X,[Y, Z]\bigr] + \bigl[ Y,[Z, X]\bigr] + \bigl[ Z,[X, Y]\bigr] = 0$$
for odd elements $X,Y,Z\in \g_{1}$ (see \cite[\S 1.6.1]{SoS} for a detailed discussion).

\label{representation-rings}

There are several different ways to extend the notion of representation rings
to Lie superalgebras. If $\g$ is a Lie superalgebra, then it is natural to
define a $\Z_{2}$-equivariant version of the representation ring constructed from $\g$-super\-modules. One complication is that, unlike finite dimensional modules over a finite dimensional complex semisimple Lie algebra, supermodules over Lie superalgebras do not in general possess complete reducibility. In order
to circumvent this problem in the representation ring, where we would normally
work with direct sums $C \cong A \oplus B$ of modules, we instead consider short
exact sequences of supermodules
\begin{equation}\label{eq:extension}
	0 \to A \to C \to B \to 0
\end{equation}
with even maps.

\begin{definition}
	The graded representation ring of a Lie superalgebra $\g$ is
	\begin{equation} \label{eq:RZ2}
		R_{\Z_{2}}(\g) :=
		F_{\Z_{2}}(\g)\,/\,
		\bigl(\text{ if $A \subset C$, then $[C/A] = [C] - [A]$ } \bigr),
	\end{equation}
	where $F_{\Z_{2}}(\g)$ is the free abelian group
	generated by (even) isomorphism classes of $\Z_{2}$-graded
	finite dimensional $\g$-modules.
\end{definition}

Recall from Section \ref{subsection-representation-ring} that a $\Z_{2}$-graded $\g$-module, or $\g$-supermodule, is a super vector space $V = V_{0} \oplus V_{1}$ with an action of $\g$ such that the action of the even component $\g_{0}$ preserves the grading, while the action of the odd component $\g_{1}$ interchanges $V_{0}$ and $V_{1}$.
By the Jordan-H\"older decomposition theorem, for any $\g$-supermodule $V$, its class $[V]$ in the Grothendieck group $R_{\Z_{2}}(\g)$ can be written uniquely
as a sum $[V] = \sum [V_{i}]$ of classes of irreducible $\g$-supermodules. Thus, the graded representation ring $R_{\Z_{2}}(\g)$
is additively the free abelian group generated by the (even) isomorphism classes of irreducible $\g$-supermodules.

The ring structure on $R_{\Z_{2}}(\g)$ is given by the tensor product, and the $\g$-action on $V\otimes W$ is
\begin{equation}\label{eq:tensor-product-action}
	X(v\otimes w) = (Xv) \otimes w + (-1)^{|X|\,|v|}\,v \otimes (Xw)
\end{equation}
for homogeneous elements $v\in V$, $w\in W$, and $X\in\g$. 
The grading operator is
\begin{equation}\label{eq:tensor-product-grading}
	(-1)^{F}_{V\otimes W} = (-1)^{F}_{V} \otimes (-1)^{F}_{W},
\end{equation}
which gives us the decomposition
\begin{equation}\label{eq:tensor-grading}
	(V\otimes W)_{0} = (V_{0} \otimes W_{0}) \oplus (V_{1} \otimes W_{1}),
	\quad
	(V\otimes W)_{1} = (V_{1} \otimes W_{0}) \oplus (V_{0} \otimes W_{1}),
\end{equation}
and the multiplicative identity is the trivial representation of $\g$
on the field $k$, viewed as a one-dimensional super vector space with no odd component.  The tensor product in fact gives a commutative product on $R_{\Z_{2}}(\g)$, since $V \otimes W$ and $W \otimes V$ are isomorphic as
$\g$-supermodules via the map $v \otimes w \mapsto (-1)^{|v|\,|w|}\,w \otimes v$.

If $\g_{0}$ is a Lie algebra with no odd component, then $\Z_{2}$-graded $\g_{0}$-modules are simply pairs of $\g_{0}$-modules, and the graded representation ring is two copies of the standard representation ring, $R_{\Z_{2}}(\g_{0}) \cong R(\g_{0}) \oplus R(\g_{0})$. On the other hand, we can also construct a representation group built out of \emph{ungraded} representations of a Lie superalgebra:

\begin{definition}
The ungraded representation group of a Lie superalgebra $\g$ is
\begin{equation}\label{eq:R0}
	R_{0}(\g) :=
	F_{0}(\g)\,/\,
	\bigl(\text{ if $A \subset C$, then $[C/A] = [C] - [A]$ } \bigr),
\end{equation}
where $F_{0}(\g)$ is the free abelian group generated by isomorphism classes of
\emph{ungraded} finite-dimen\-sional $\g$-modules.
\end{definition}

The ungraded representation group is a module over the graded representation ring, with product
\begin{equation}\label{eq:module}
\otimes : R_{\Z_{2}}(\g) \times R_{0}(\g) \to R_{0}(\g)
\end{equation}
induced by the tensor product. If $V$ is a $\g$-supermodule and $W$ a $\g$-module, 
then the $\g$-action on $V \otimes W$ is once again given by (\ref{eq:tensor-product-action}).  We note that (\ref{eq:tensor-product-action}) requires a $\Z_{2}$-grading only on the first factor $V$.  On the other hand, the grading operator (\ref{eq:tensor-product-grading}) requires $\Z_{2}$-gradings on both factors $V$ and $W$,
so here the tensor product $V \otimes W$ is a $\g$-module, not a $\g$-supermodule.

Although we can take tensor products of two ungraded $\g$-modules as vector spaces, if the odd component $\g_{1}$ acts nontrivially on the factors, then the tensor product does not necessarily admit a $\g$-action. Thus, we do not have a ring structure in the ungraded case.  If $\g_{0}$ is a Lie algebra with no odd component, then this definition  gives the usual representation ring $R_{0}(\g_{0}) = R(\g_{0})$. When working with Lie superalgebras, it is uncommon to consider ungraded representations. However, we can express the ungraded representation group $R_{0}(\g)$ in terms of graded representations, as we will show in Section \ref{redux}.

By comparing the ungraded representation group and the graded representation ring, we construct another variant of the representation ring, which is closely related to equivariant $K$-theory. There is a standard functor $f: R_{\Z_{2}}(\g) \to R_{0}(\g)$, taking $\g$-super\-modules to $\g$-modules by simply forgetting their gradings. However, for now we are more interested in the adjoint functor $\Delta: R_{0}(\g) \to R_{\Z_{2}}(\g)$ going in the opposite direction. Given a representation $r : \g \to \End(U)$ of the Lie superalgebra $\g$ on an ungraded vector space $U$, we construct a $\Z_{2}$-graded representation $$\Delta r : \g \to \End(U\oplus U)$$ on the super vector space $\Delta U = U \oplus  U$ (i.e., with $(\Delta U)_{0} = (\Delta U)_{1} = U$), by
\begin{equation}\label{eq:deltar}
(\Delta r) (X_{0}) := 
	\begin{pmatrix} r(X_{0}) & 0 \\ 0 & r(X_{0}) \end{pmatrix} 
	\quad\text{and}\quad
	(\Delta r) (X_{1}) := 
	\begin{pmatrix} 0 & r(X_{1}) \\ r(X_{1}) & 0 \end{pmatrix},
\end{equation}
for $X_{0}\in \g_{0}$ and $X_{1}\in \g_{1}$.
If $\g_{0}$ is a Lie algebra, then
$$\Delta : \bigl( R_{0}(\g_{0}) \cong R(\g_{0}) \bigr) \longrightarrow \bigl( R_{\Z_{2}}(\g_{0}) \cong R(\g_{0}) \oplus R(\g_{0}) \bigr)$$
is the diagonal map, hence our choice of notation.  For a Lie superalgebra with a nontrivial odd component, the homomorphism $\Delta$ is the diagonal map only at the level of vector spaces, but not at the level of supermodules.

\begin{definition}The super representation ring of a Lie superalgebra
$\g$ is the quotient
\begin{equation}\label{eq:SR}
	SR(\g) := R_{\Z_{2}}(\g) / \Delta R_{0}(\g)
\end{equation}
of the ring of $\g$-supermodules by the image of the ring of ungraded $\g$-modules under the diagonal map $\Delta$.
\end{definition}

For a Lie algebra $\g_{0}$ with no odd part, the super representation ring is once again isomorphic to the usual representation ring: $SR(\g_{0}) \cong R(\g_{0})$. We now investigate the product structure on the super representation ring.

\begin{lemma}\label{tensor}
	If $V$ is a $\g$-supermodule and $U$ is an ungraded
	$\g$-module, then
	$$\Delta ( V \otimes U ) \cong V \otimes \Delta U.
	$$
\end{lemma}
\begin{proof}
As vector spaces, we have $(\Delta U)_{0} = (\Delta U)_{1} = U$. It follows
that the even and odd components of $\Delta(V\otimes U)$ are precisely equal
to the even and odd components of $V \otimes \Delta U$ as vector spaces:
\begin{align*}
\bigl(\Delta(V\otimes U)\bigr)_{0} &= V \otimes U  = ( V_{0} \otimes U) \oplus ( V_{1} \otimes U) \\
&= \bigl(V_{0} \otimes (\Delta U)_{0}\bigr)
	\oplus \bigl(V_{1} \otimes (\Delta U)_{1}\bigr)
= (V \otimes \Delta U)_{0}, \\
\bigl(\Delta(V\otimes U)\bigr)_{1} &= V \otimes U = (V_{0} \otimes U) \oplus (V_{1} \otimes U) \\
&= \bigl(V_{0} \otimes (\Delta U)_{1}\bigr)
	\oplus \bigl(V_{1} \otimes (\Delta U)_{0}\bigr)
= (V \otimes \Delta U)_{1}.
\end{align*}
Consider the map $\Delta(V\otimes U)\to V \otimes \Delta U$ which
acts as the identity on these equivalent super vector spaces. This map
clearly commutes with the action of $\g$ on the tensor product,
and so it is a $\g$-supermodule isomorphism.
\end{proof}

\begin{corollary}The image $\Delta R_{0}(\g) \subset R_{\Z_{2}}(\g)$ of the diagonal map is an ideal, and thus the quotient
$SR(\g) = R_{\Z_{2}}(\g) / \Delta R_{0}(\g)$ is indeed a ring.
\end{corollary}

Our definition (\ref{eq:SR}) of the super representation ring is motivated by the following analogue of Proposition \ref{prop-1} (see Section~\ref{parity-reversal} for our definition of the right parity reversal $V^{\Pi}$, which has the opposite $\Z_{2}$-grading but the same underlying $\g$-action as $V$, without introducing any signs):

\begin{proposition}\label{prop-4}
The image $\Delta R_{0}(\g)$ of the diagonal map is generated by the classes of $\g$-super\-modules $[V]\in R_{\Z_{2}}(\g)$ which admit a $\g$-equivariant parity reversing involution $\alpha: V\to V^{\Pi}$.
\end{proposition}

\begin{proof}
From the description (\ref{eq:deltar}) of $\Delta U$, we have an obvious choice for the $\g$-equivariant parity reversing involution on $\Delta U$, namely
$$\alpha = \begin{pmatrix} 0 & \mathrm{Id} \\ \mathrm{Id} & 0\end{pmatrix}.$$
Conversely, suppose we are given a representation $r : \g \to \End(V)$ of $\g$ on a super vector space $V = V_{0} \oplus V_{1}$, equipped with a $\g$-equivariant involution $\alpha$ which interchanges $V_{0}$ and $V_{1}$. We can now construct an ungraded $\g$-action on the even component $p_{0}(V) = V_{0}$, given by
\begin{equation}\label{eq:p0}
	(p_{0}r)(X_{0}) = \bigl[r(X_{0}) \bigr]_{V_{0}} \quad\text{and}\quad
	(p_{0}r)(X_{1}) = \bigl[\alpha \, r(X_{1}) \bigr]_{V_{0}}
\end{equation}
for $X_{0}\in \g_{0}$ and $X_{1}\in\g_{1}$. (We could likewise define an ungraded $\g$-action on the odd component $p_{1}(V) = V_{1}$,
and the involution $\alpha$ then gives us an isomorphism $p_{0}(V) \cong p_{1}(V)$ as ungraded $\g$-modules.) Now, applying the diagonal map, we have $\Delta(p_{0}V) = V_{0} \oplus V_{0}$, which is isomorphic to $V = V_{0} \oplus V_{1}$ via the map $1 \oplus \alpha$. Combining (\ref{eq:deltar}) and (\ref{eq:p0}), the $\g$-action on
$\Delta(p_{0}V)$ is then
$$\bigl(\Delta (p_{0}r)\bigr) (X_{0}+X_{1}) = 
	\begin{pmatrix} \bigl[r(X_{0})\bigr]_{V_{0}} & \bigl[\alpha\,r(X_{1})\bigr]_{V_{0}} \\ \bigl[\alpha\,r(X_{1})\bigr]_{V_{0}} & \bigl[r(X_{0})\bigr]_{V_{0}} \end{pmatrix} 
$$
for $X_{0}\in \g_{0}$ and $X_{1}\in \g_{1}$, and since $\alpha$ commutes with the $\g$-action, we have
$$\begin{pmatrix} 1 & 0 \\ 0 & \alpha \end{pmatrix}
  \bigl(\Delta (p_{0}r)\bigr)(X)
  \begin{pmatrix} 1 & 0 \\ 0 & \alpha \end{pmatrix}
  = r(X)$$
for $X\in \g$. Thus $\Delta (p_{0} V) \cong V$, and so the class $[V]$ of any such $\g$-supermodule with a parity reversing involution $\alpha$ is in the image of the diagonal map $\Delta$.
\end{proof}

A comparison of Proposition \ref{prop-4} with the $n=0$ case of Proposition
\ref{prop-1} indicates that there is a Clifford algebraic interpretation of the super representation ring, which we will investigate further in Sections \ref{degree-shift} and
\ref{redux}.
We can reverse our construction of the super representation ring, substituting the forgetful map $f : R_{\Z_{2}}(\g) \to R_{0}(\g)$ for the diagonal map $\Delta : R_{0}(\g) \to R_{\Z_{2}}(\g)$. In analogy to $K$-theory, the resulting quotient gives a group which behaves much like a degree-shifted version of the super representation ring $SR(\g)$.

\begin{definition}The ungraded super representation group of a Lie superalgebra
$\g$ is the quotient
\begin{equation}\label{eq:SR1}
	SR'(\g) := R_{0}(\g) / f\, R_{\Z_{2}}(\g)
\end{equation}
of the group of ungraded $\g$-modules by the image of the ring of $\g$-supermodules under the forgetful map $f$.
\end{definition}

For any Lie algebra $\g_{0}$, the forgetful map $f$ is surjective, and so the ungraded super representation group $SR'(\g_{0}) = 0$ vanishes. 

For the remainder of this section we work over the complex numbers.
We recall from the above discussion that the tensor product of a supermodule with an ungraded module is once again an ungraded $\g$-module with $\g$-action given by (\ref{eq:tensor-product-action}). The tensor product then gives 
the ungraded representation group $R_{0}(\g)$ the structure (\ref{eq:module})
of a module over the graded representation ring $R_{\Z_{2}}(\g)$. To show that this product structure descends to the ungraded super representation group $SR'(\g)$, we require the following technical lemma:

\begin{lemma}\label{lemma:product}
	Given two complex ungraded $\g$-modules $U$ and $W$, the ungraded tensor product $\g$-module $(\Delta U) \otimes W$ nevertheless admits a $\Z_{2}$-grading.
\end{lemma}

\begin{proof}
As vector spaces, we have $(\Delta U)\otimes W = (U \otimes W) \oplus (U \otimes W)$, and the $\g$-action is given by
\begin{align*}
X_{0} &\mapsto \begin{pmatrix} X_{0} \otimes 1 + 1 \otimes X_{0} & 0 \\ 0 & X_{0} \otimes 1 + 1 \otimes X_{0}\end{pmatrix}, \\
X_{1} &\mapsto \begin{pmatrix} 1 \otimes X_{1} & X_{1} \otimes 1 \\ X_{1} \otimes 1 & -1 \otimes X_{1}\end{pmatrix},
\end{align*}
for $X_{0}\in \g_{0}$ and $X_{1}\in \g_{1}$. We can then consider the grading operator
\begin{equation}\label{eq:new-grading}
(-1)^{F} = \begin{pmatrix} 0 & i\,\mathrm{Id} \\ -i\,\mathrm{Id} & 0 \end{pmatrix},
\end{equation}
which commutes with the action of $\g_{0}$ and anti-commutes with the action of $\g_{1}$. 
\end{proof}

Given two ungraded $\g$-modules $U$ and $W$, we can define their product to be the $\g$-supermodule
$$U \boxtimes W := (\Delta U) \otimes W$$
with the $\Z_{2}$-grading given by (\ref{eq:new-grading}). This product induces
a map
\begin{equation}\label{eq:twisted-product}
\boxtimes : R_{0}(\g) \times R_{0}(\g) \to R_{\Z_{2}}(\g),
\end{equation}
which we now show descends to the super representation ring.

\begin{proposition}
Over the complex numbers, the direct sum $SR(\g) \oplus SR'(\g)$ forms a supercommutative $\Z_{2}$-graded ring with $SR(\g)$ in degree 0 and $SR'(\g)$ in degree 1. In particular, we have:
\begin{enumerate}
\item The ungraded super representation group $SR'(\g)$ is a module over the super representation ring $SR(\g)$ with product $\otimes : SR(\g) \times SR'(\g) \to SR'(\g)$ induced by the tensor product (\ref{eq:module}).
\item The product (\ref{eq:twisted-product}) induces a product $\boxtimes: SR'(\g) \times SR'(\g) \to SR(\g)$.
\item The product of two ungraded $\g$-modules $U$ and $W$ satisfies the identity$$U \boxtimes W \cong (W \boxtimes U)^{\Pi},$$ and thus $[U
] \boxtimes [W] = - [W] \boxtimes [U]$ is anti-commutative in the super representation ring.
\end{enumerate}
\end{proposition}

\begin{proof}
We begin with (1). If $V$ and $W$ are $\g$-supermodules, then we clearly have $V \otimes f(W) = f(V\otimes W)$ as $\g$-modules. Indeed, the $\g$-actions on these two tensor products are both given by (\ref{eq:tensor-product-action}). On the other hand,
if $U$ and $W$ are ungraded $\g$-modules, then by Lemma~\ref{lemma:product} we see that the ungraded $\g$-module $(\Delta U) \otimes W$ admits a $\Z_{2}$-grading and is therefore in the image of the forgetful map $f$. In summary, the tensor product satisfies
\begin{align*}
R_{\Z_{2}}(\g) \otimes f R_{\Z_{2}}(\g) &\subset f R_{\Z_{2}}(\g), \\
\Delta R_{0}(\g) \otimes R_{0}(\g) &\subset f R_{\Z_{2}}(\g),
\end{align*}
and it therefore descends to a product $SR(\g) \times SR'(\g) \to SR'(\g)$.

We will only sketch the proofs of (2) and (3) here, as we prove the same results in the language of Clifford algebraic degree shifts in Proposition~\ref{SR-tensor} below. We establish the equivalence of these two approaches in Corollary~\ref{equivalent}. For (3), we note that the underlying vector spaces for $U \boxtimes W$ and $W \boxtimes U$ are isomorphic by interchanging the two tensor factors, and that reversing the order corresponds to choosing the opposite sign for the grading operator (\ref{eq:new-grading}). For (2), in order that the product (\ref{eq:twisted-product}) descend to the super representation ring, we must show that for an ungraded $\g$-module $U$ and a $\g$-supermodule $V$, the product $U \boxtimes fV$ lies in the image of $\Delta$. (The analogous result for $fV \boxtimes U$ then follows from the anti-commutativity of the product.) However,
$$(\Delta U) \otimes fV
\cong f \bigl( (\Delta U) \times V \bigr)
\cong f \bigl( V \times \Delta U )
\cong f \Delta ( V \otimes U ),$$
using the commutativity of the tensor product of $\g$-supermodules and Lemma~\ref{tensor}. Removing the forgetful map $f$ on the right hand side and reestablishing the $\Z_{2}$-grading on $(\Delta U)\otimes fV$, we can then show that $U \boxtimes fV \cong \Delta( V \otimes U)$, and thus $R_{0}(\g) \boxtimes f R_{\Z_{2}}(\g) \subset \Delta R_{0}(\g)$.
\end{proof}

%

\section{Involutions}
\label{involutions}

In order to study the diagonal map $\Delta : R_{0}(\g) \to R_{\Z_{2}}(\g)$, the forgetful map $f : R_{\Z_{2}}(\g) \to R_{0}(\g)$, the super representation ring $SR(\g)$, and the ungraded super representation group $SR'(\g)$ in more detail, we introduce involutions on both the graded representation ring and the ungraded representation group. In this section we continue working over an arbitrary field of characteristic 0. Readers who are experienced with Lie superalgebras and interested in the complex case, in which Schur's Lemma simplifies matters significantly, are directed to skip directly to Section \ref{degree-shift}, in which we introduce Clifford algebraic degree shifts. What follows is intended for readers new to Lie superalgebras, readers interested in the real, non-algebraically closed, case, or readers interested in ungraded representations.  

\subsection{Parity reversal}
\label{parity-reversal}

Let $\Pi$ be the parity reversal operator on super vector spaces, which leaves the underlying vector space unchanged but switches the $\Z_{2}$-grading: $(V^{\Pi})_{0} = V_{1}$ and $(V^{\Pi})_{1} = V_{0}$. Let $\pi : V \to V^{\Pi}$ be the canonical odd homomorphism leaving the underlying vectors fixed.
If $V$ is a $\g$-supermodule, then we define a $\g$-supermodule structure on $V^{\Pi}$ by taking
\begin{equation}\label{eq:right-pi}
	X \cdot (v^{\pi}) := (X \cdot v)^{\pi}
\end{equation}
for all $X\in \g$ and $v\in V$. Equivalently, we have $V^{\Pi} := V \otimes_{\Z} \Pi\Z.$

\begin{remark}\label{remark:pi}
It is important to note that our convention (\ref{eq:right-pi}) differs from that of the more standard left parity reversal operator $\Pi(V)$, which satisfies $X \cdot \pi(v) = (-1)^{|X|}\,\pi(X\cdot v)$. Rather, we are using the right parity reversal operator, sometimes denoted $(V)\Pi$, which commutes (in the usual, non-supersymmetric sense) with the left $\g$-action. See \cite[\S 1.2.2]{SoS}. This choice of sign convention does not affect the representation ring, as the left and right parity reversals are conjugate to each other in the sense of Section~\ref{section:conjugation} below, and thus $\Pi(V) \cong (V)\Pi$ by Remark~\ref{remark:conjugation}. On the other hand, it is vital to recall our sign conventions when we refer to a parity reversing involution $\alpha: V \to V^{\Pi}$. Our choice of the right parity reversal operator agrees with the ${}^{*}$ operator in \cite{ABS} and allows us to reproduce the $K$- and $KO$-theory of a point as the super representation rings of Clifford algebras.
\end{remark}

The operator $\Pi$ gives an involution on $F_{\Z_{2}}(\g)$ which descends to an involution on the graded representation ring $R_{\Z_{2}}(\g)$.

\begin{definition}
	The self-dual $\Rp(\g)$ and anti-dual $\Rm(\g)$ representation 
	ideals
	of a Lie superalgebra $\g$ are the $+1$ and $-1$ eigenspaces of the parity
	reversal operator $\Pi:R_{\Z_{2}}(\g)\to R_{\Z_{2}}(\g)$.
\end{definition}

Alternatively, we could define $\Rpm(\g)$ in terms of the $\pm 1$ eigenspaces
of the parity reversal operator $\Pi: F_{\Z_{2}}(\g) \to F_{\Z_{2}}(\g)$ acting on the free abelian group generated by isomorphism classes of $\g$-supermodules,
using the same construction as in (\ref{eq:RZ2}). The self-dual and anti-dual representation ideals both agree with the standard representation ring,
$$\Rp(\g_{0}) \cong \Rm(\g_{0}) \cong R(\g_{0}),$$
for a Lie algebra $\g_{0}$.

\begin{proposition}\label{ideal}
	The groups $\Rp(\g)$ and $\Rm(\g)$ are indeed ideals in the graded representation
	ring $R_{\Z_{2}}(\g)$.
\end{proposition}
\begin{proof}
The key point here is that if we take the parity reversal of the tensor product $V \otimes W$ of two $\g$-supermodules $V$ and $W$, it follows from the decomposition (\ref{eq:tensor-grading}) that
\begin{equation}\label{eq:pi-tensor}
( V\otimes W)^{\Pi} \cong V^{\Pi} \otimes W \cong V \otimes W^{\Pi}.
\end{equation}
From this, it is clear that the $\pm 1$ eigenspaces of $\Pi$ are closed under the tensor products, and thus $\Rpm(\g)$ are indeed (non-unital) subrings of $R_{\Z_{2}}(\g)$.
If $V \cong V^{\Pi}$ is a self-dual $\g$-supermodule corresponding to a class $[V]\in \Rp(\g)$, then for any $\g$-supermodule $W$ we have
$( V \otimes W )^{\Pi} \cong V^{\Pi} \otimes W \cong V \otimes W$, and thus
the class of the product $[V][W] = [V \otimes W]$ also lies in $\Rp(\g)$.
Similarly, if $[V] = -\Pi[V] \in \Rm(\g)$ is a virtual anti-dual $\g$-supermodule, then in light of (\ref{eq:pi-tensor}),
the product $[V][W]$ for any virtual $\g$-supermodule $[W] \in R_{\Z_{2}}(\g)$ is likewise anti-dual, since
$\Pi([V][W]) = (\Pi[V])[W] = -[V][W]$.
It follows that both $\Rp(\g)$ and $\Rm(\g)$ are ideals in $R_{\Z_{2}}(\g)$. 
\end{proof}

If $V$ is an irreducible $\g$-supermodule, corresponding to a generator of the graded representation ring $R_{\Z_{2}}(\g)$, then the parity reversed $\g$-supermodule $V^{\Pi}$ is likewise irreducible.  (If $V^{\Pi}$ were reducible, then applying $\Pi$ to a short exact sequence (\ref{eq:extension}) decomposing $V^{\Pi}$, we would obtain a similar short exact sequence for $V$,
which contradicts the irreducibility of $V$.) There are now two possibilities: either $V$ is isomorphic to its parity reversal, or its parity reversal lies
in a distinct isomorphism class. (These two types of irreducible supermodules were first identified by Bernstein and Leites in \cite{BL} and are commonly referred to as 
type Q (for queer) and type M (for matrix), respectively. In the complex case, this classification of irreducibles is related to Schur's Lemma, where the type M and type Q irreducibles have 1-dimensional and $1|1$-dimensional endomorphism algebras, respectively.) In these terms, the self-dual representation
ideal $\Rp(\g)$ is generated by the classes of (type Q) irreducible $\g$-supermodules $V$ with $V^{\Pi} \cong V$, as well as by the classes of $\g$-supermodules of the form $V \oplus V^{\Pi}$
for (type M) irreducible $V$ with $V^{\Pi} \ncong V$.  On the other hand,
the anti-dual representation ideal $\Rm(\g)$ is generated by virtual $\g$-supermodules of the form $[V] - [V^{\Pi}]$ for (type M) irreducible $V$ with $V^{\Pi} \ncong V$.

Recall from the introduction the 
version of the representation ring (see \cite{BK0}) built from isomorphism classes of finite dimensional  supermodules where we identify a supermodule with its parity reversal: $[V] \sim [V^{\Pi}]$. We can construct this representation ring
by taking the quotient of the graded representation ring $R_{\Z_{2}}(\g)$ by
the anti-dual representation ideal $\Rm(\g)$, generated by the differences $[V]-[V^{\Pi}]$.  (Actually, the product structure on this representation ring is slightly different, with the product of two type Q irreducible supermodules being half their tensor product. See \cite{ShM,Se2} for a complete discussion.)  Equivalently, we can construct this representation ring using
the method of (\ref{eq:RZ2}), but relaxing the condition on $F_{\Z_{2}}(\g)$ that all isomorphisms must be even maps. We then have:

\begin{proposition}The representation ring of $\g$-supermodules up to parity reversal is the quotient
$$R_{\Z_{2}}(\g) / \Rm(\g) \cong F(\g) /
\bigl( \text{ if $A \subset C$, then $[C/A] = [C] - [A]$ } \bigr),$$
where $F(\g)$ is the free abelian group generated by
isomorphism classes of $\g$-super\-modules with respect to both even and
odd isomorphisms.
\end{proposition}

In light of Proposition \ref{ideal}, the representation ring $R_{\Z_{2}}(\g) / \Rm(\g)$ is indeed a ring with respect to the tensor product. On the other hand, if we consider the quotient of the graded representation ring by the self-dual representation ideal, we again obtain a ring, and we have:

\begin{lemma}\label{anti-sym}
The anti-symmetrizing map $1-\Pi: R_{\Z_{2}}(\g) \to \Rm(\g)$, given by $1-\Pi:[V] \mapsto [V] - [V^{\Pi}]$, is surjective and descends to an isomorphism
\begin{equation*}
	1-\Pi: R_{\Z_{2}}(\g) / \Rp(\g) \xrightarrow{\cong} \Rm(\g)
\end{equation*}
as additive groups. (This map is not a ring homomorphism.)
\end{lemma}

\begin{proof}
The generators $[V] - [V^{\Pi}]$ of $\Rm(\g)$ are clearly all in the image
of $1-\Pi$, and the kernel is generated by $[V]$ such that $[V^{\Pi}] = [V]$. 
\end{proof}

The anti-symmetrizing map $1-\Pi$ is essentially multiplication by 2 on the quotient $R_{\Z_{2}}(\g)/R_{+}(\g)$, because the parity reversal operator $\Pi$ acts by $-1$ there.  In contrast, the analogous symmetrizing
map $R_{\Z_{2}}(\g) \to \Rp(\g)$ induced by the map $V \mapsto V \oplus V^{\Pi}$ on $\g$-supermodules is \emph{not} surjective.  In particular, if
$V$ is an irreducible $\g$-super\-module with $V\cong V^{\Pi}$ via an (even) isomorphism,
thereby giving a generator of the self-dual representation ideal $\Rp(\g)$, then we cannot decompose $V$ as $W \oplus W^{\Pi}$. Nevertheless, one could indeed construct a group (but not ring) isomorphism $$R_{\Z_{2}}(\g)/\Rm(\g) \to \Rp(\g)$$ by considering the generators $[V]$ for $\g$-supermodules with $V\cong V^{\Pi}$ and $V\ncong V^{\Pi}$
separately. However, such a map would not be canonical, and we will not consider it further.

\subsection{Conjugation}\label{section:conjugation}

We next introduce an involution on the ungraded representation group $R_{0}(\g)$ which is similar to parity reversal for graded representations. If $r : \g \to \End(U)$ is an ungraded representation of $\g$ on a vector space $U$, we construct a new ungraded representation of $\g$ on $U$ by taking
\begin{equation}\label{eq:rdag}
r^{\dag}(X_{0}) = r(X_{0}), \qquad 
  r^{\dag}(X_{1}) = - r(X_{1}),
\end{equation}
for $X_{0}\in\g_{0}$ and $X_{1}\in\g_{1}$. In other words, if $(-1)^{F}$
is the grading operator on the Lie superalgebra $\g$, then we put
$r^{\dag}(X) = r\bigl( (-1)^{F} X \bigr)$ for all $X\in \g$. For notational
purposes, we let $U^{\dag}$ denote the new ungraded $\g$-module, which has the same
underlying vector space as $U$ but with the new action $r^{\dag}$. This map $U\mapsto U^{\dag}$ gives an involution of the ungraded representation group $R_{0}(\g)$ and is analogous to complex conjugation. (Indeed, if we extend the definition (\ref{eq:rdag}) to representations of associative superalgebras and consider ungraded real representations of the Clifford algebra $\mathrm{Cl}(1)\cong \C$, then this map takes a complex vector space to the complex conjugate vector space.) Just as we defined the self-dual and anti-dual representation ideals, we can do the same for the conjugation operator:

\begin{definition}
The self-conjugate $\Rsc(\g)$ and the anti-conjugate $\Rac(\g)$ representation groups of a Lie superalgebra $\g$ are the $+1$ and $-1$ eigen\-spaces, respectively, of the conjugation operator $\dag : R_{0}(\g) \to R_{0}(\g)$.
\end{definition}

Alternatively, we could define $\Rsc(\g)$ and $\Rac(\g)$ in terms of
the $+1$ and $-1$ eigenspaces of the conjugation operator $\dag:F_{0}(\g) \to F_{0}(\g)$ acting on the free abelian group generated by isomorphism classes of ungraded $\g$-modules, using the same construction as in (\ref{eq:R0}).
For a Lie algebra $\g_{0}$ with no odd component, all $\g_{0}$-modules are self-conjugate, and thus $$\Rsc(\g_{0}) = R_{0}(\g_{0}) = R(\g_{0}),$$ while $\Rac(\g_{0}) = 0.$

As we did with the self-dual and anti-dual representation ideals, we can describe the generators of the self-conjugate and anti-conjugate representation groups in terms of irreducibles. First, note
that if $U$ is an irreducible ungraded $\g$-module, then its conjugate
$U^{\dag}$ is likewise irreducible. The self-conjugate representation group $\Rsc(\g)$ is then generated by the classes of irreducible self-conjugate $\g$-modules $U\cong U^{\dag}$, as well as the classes of
$\g$-modules of
the form $U \oplus U^{\dag}$, for irreducible $\g$-modules $U$ such
that $U\ncong U^{\dag}$.  On the other hand, the anti-conjugate representation group $\Rac(\g)$ is generated by virtual representations
of the form $[U] - [U^{\dag}]$ for irreducible $\g$ modules $U$ with $U\ncong U^{\dag}$. Lemma \ref{anti-sym} also holds with conjugation in place of parity reversal:

\begin{lemma}\label{anti-sym-dag}
The anti-symmetrizing map $1-\dag: R_{0}(\g) \to \Rac(\g)$, given by $$1-\dag:[U] \mapsto [U] - [U^{\dag}],$$ is surjective and descends to a group isomorphism
\begin{equation*}
	1-\dag: R_{0}(\g) / \Rsc(\g) \xrightarrow{\cong} \Rac(\g).
\end{equation*}
\end{lemma}

\begin{remark}\label{remark:conjugation}
We observe that the conjugation operator can also be applied to $\g$-super\-modules, although
it acts trivially on the graded representation ring $R_{\Z_{2}}(\g)$. 
This is because a $\g$-supermodule admits a $\Z_{2}$-action which intertwines with the $\Z_{2}$-action on the Lie superalgebra $\g$, and the corresponding grading operator gives isomorphisms between $\g$-supermodules and their conjugates. More precisely,
given a $\Z_{2}$-graded representation $r : \g \to \End(V)$ of $\g$ on a super vector space $V$, then we have
\begin{equation*}
	r^{\dag}(X) = r\bigl( (-1)^{F}X \bigr) = (-1)^{F} \circ r(X) \circ (-1)^{F}
\end{equation*}
for all $X\in \g$, and therefore $V \cong V^{\dag}$. In analogy to Proposition~\ref{prop-4}, we observe that the image of the forgetful map consists of those ungraded $\g$-modules $U$ which admit conjugate $\g$-equivariant involutions $\beta : U \to U^{\dag}$.
\end{remark}

\section{The periodic exact sequence}
\label{les}

In this section, we consider only additive group structures, and all maps here are group homomorphisms.  If we ignore products, then the super representation ring $SR(\g)$ and the ungraded super representation group $SR'(\g)$ constructed in Section \ref{representation-rings} are essentially mirror images of one another.  In particular, every statement we make in this section about
the diagonal map $\Delta : R_{0}(\g) \to R_{\Z_{2}}(\g)$ corresponds to a statement about the forgetful map $f : R_{\Z_{2}}(\g) \to R_{0}(\g)$, but with the involutions $\Pi$ and $\dag$ interchanged (recall from Remark~\ref{remark:pi} that we are using the right parity reversal operator, which maintains the same underlying module structure, without introducing any signs). It is only once we consider products that the ring $SR(\g)$ exhibits different behavior from its ungraded version $SR'(\g)$.  We also continue working over an arbitrary field of characteristic 0.

\subsection{Kernels and images}
Recall from linear algebra that the complexification of an already complex vector space is isomorphic to the direct sum of the original vector space and its complex conjugate. The analogous statement holds in the context of the conjugation map $\dag$ acting on ungraded $\g$-modules in $R_{0}(\g)$, where
the role of complexification is played by the diagonal map $\Delta$.

\begin{proposition}\label{composition}
The compositions of the diagonal map $\Delta : R_{0}(\g) \to R_{\Z_{2}}(\g)$ with the forgetful map $f: R_{\Z_{2}}(\g) \to R_{0}(\g)$, taken in either order, correspond to the symmetrizing operators $1+\Pi$ and $1+\dag$ on $R_{\Z_{2}}(\g)$ and $R_{0}(\g)$, respectively. In other words, the following diagram commutes:
$$
\begin{diagram}[vtrianglewidth=1.73em]
R_{\Z_{2}}(\g) & &\rTo^{1+\Pi} & & R_{\Z_{2}}(\g) \\
&\rdTo>{f} & &\ruTo<{\Delta} & & \rdTo>{f} \\
&& R_{0}(\g) && \rTo^{1+\dag} && R_{0}(\g)
\end{diagram}
$$
Consequently, the symmetrizing operators commute with the diagonal map and the forgetful map.
\end{proposition}
\begin{proof}
	First we consider the bottom right triangle, the composition $$R_{0}(\g)\xrightarrow{\Delta}R_{\Z_{2}}(\g) \xrightarrow{f} R_{0}(\g).$$
	Let $r : \g \to \End(U)$ be an ungraded representation of $\g$ on
	a vector space $U$. 
	Then in the graded $\g$-module
	$\Delta U = U \oplus U = U \otimes k^{2}$ ($k$ is our field),
	consider the diagonal and anti-diagonal subspaces
	$$U_{+} = \left\{ \begin{pmatrix}u \\ u\end{pmatrix} : u\in U \right\}, \qquad
	  U_{-} = \left\{ \begin{pmatrix}\phantom{-}u \\ -u\end{pmatrix} : u\in U \right\}.$$
	We clearly have the decomposition $\Delta U = U_{+} \oplus U_{-}$ as ungraded
	vector spaces.
	With respect to the
	action of
	$\Delta r$ given by (\ref{eq:deltar}), both $U_{+}$ and $U_{-}$
	are $\g$-invariant subspaces, with actions given by
	\begin{align*}
		(\Delta r) (X_{0}+X_{1})|_{U_{+}\;} &:
		\begin{pmatrix}u \\ u \end{pmatrix}
		\mapsto \begin{pmatrix}r(X_{0}+X_{1})(u) \\ r(X_{0}+X_{1})(u)\end{pmatrix}, \\
		(\Delta r) (X_{0}+X_{1})|_{U_{-}} &:
		\begin{pmatrix} \phantom{-}u \\ -u\end{pmatrix} \mapsto
		\begin{pmatrix} \phantom{-}r(X_{0}-X_{1})(u) \\ - r(X_{0}-X_{1})(u) \end{pmatrix}.
	\end{align*}
	We then see that $U_{+} \cong U$ and by (\ref{eq:rdag}) that
	$U_{-} \cong U^{\dag}$, and thus $\Delta U \cong U \oplus U^{\dag}$ as
	ungraded $\g$-modules.
	
	Next, we consider the top left triangle, the composition
	$$R_{\Z_{2}}(\g) \xrightarrow{f} R_{0}(\g)\xrightarrow{\Delta} R_{Z_{2}}(\g).$$
	Given a $\g$-supermodule $V = V_{0}\oplus V_{1}$, forgetting its
	grading and
	applying the diagonal map gives
	$$\Delta V = V \oplus V = 	(V_{00}\oplus V_{10}) \oplus ( V_{01} \oplus V_{11}),$$
	where the even part $\g_{0}$ preserves all four components, while the odd part $\g_{1}$
	interchanges the components in the pairs $V_{00} \leftrightarrow V_{11}$ and $V_{10}\leftrightarrow V_{01}$. It follows that $V_{00}\oplus V_{11}$ and $V_{10}\oplus V_{01}$ are complementary 
	$\g$-invariant
	subspaces inside $\Delta V$, which are isomorphic as $\g$-supermodules to $V$ and $V^{\Pi}$,
	 respectively. Thus $\Delta V \cong V \oplus V^{\Pi}$ in $R_{\Z_{2}}(\g)$.
\end{proof}

We now examine the diagonal map $\Delta$ and the forgetful map $f$ in more detail, considering their images and kernels.

\begin{proposition}\label{image}
The images of the diagonal map $\Delta : R_{0}(\g) \to R_{\Z_{2}}(\g)$ and the forgetful map $f: R_{\Z_{2}}(\g) \to R_{0}(\g)$ are contained in
the self-dual representation ideal and the self-conjugate representation group of $\g$, respectively:
$$\mathrm{Im}\,\Delta \subset \Rp(\g) \subset R_{\Z_{2}}(\g)
\quad \text{and} \quad
\mathrm{Im}\,f \subset \Rsc(\g) \subset R_{0}(\g).$$
\end{proposition}

\begin{proof}
Given an ungraded
$\g$-module $U$, consider the operator
$$\alpha = \begin{pmatrix} 0 & \mathrm{Id} \\ \mathrm{Id} & 0 \end{pmatrix}$$
acting on the $\g$-supermodule $\Delta U = U\oplus U$. This map $\alpha$ clearly interchanges the odd and even components of $\Delta U$, and comparing with (\ref{eq:deltar}),
we see that $\alpha$ commutes with the actions of both the odd and even components $\g_{0}$ and $\g_{1}$ of the Lie superalgebra $\g$.  The map $\alpha$ therefore gives an even $\g$-supermodule isomorphism $\Delta U \cong (\Delta U)^{\Pi}$, and thus
$[\Delta U] \in \Rp(\g)$.  As for the the forgetful map, we recall from Remark~\ref{remark:conjugation} that every $\g$-supermodule is automatically self-conjugate.
\end{proof}

\begin{proposition}\label{kernel}
The kernels of the diagonal map $\Delta : R_{0}(\g) \to R_{\Z_{2}}(\g)$ and the forgetful map $f:R_{\Z_{2}}(\g) \to R_{0}(\g)$
are
$$
\mathrm{Ker}\, \Delta = \Rac(\g) \subset R_{0}(\g)
\quad \text{and} \quad
\mathrm{Ker}\, f = \Rm(\g) \subset R_{\Z_{2}}(\g),
$$
the anti-conjugate representation group and anti-dual representation ideal of $\g$,
respectively.
\end{proposition}

\begin{proof}
Let $r : \g \to \End(U)$ be an ungraded representation of $\g$. Applying the diagonal map $\Delta : R_{0}(\g) \to R_{\Z_{2}}(\g)$,
the grading operator
$$(-1)^{F} = \begin{pmatrix} \mathrm{Id} & 0 \\ 0 & \mathrm{-Id} \end{pmatrix}$$
on $\Delta U$ and $\Delta U^{\dag}$ intertwines with the $\g$-actions of $\Delta r$ and $\Delta r^{\dag}$ given by (\ref{eq:deltar}) and
thus gives an (even) isomorphism $\Delta U \cong \Delta U^{\dag}.$ We can also
see this by noting that conjugation commutes with the diagonal map: $\Delta (U^{\dag}) = (\Delta U)^{\dag}$.
However, as $\Delta U$ is a $\g$-supermodule, we recall from Remark~\ref{remark:conjugation} that $(\Delta U)^{\dag} \cong \Delta U$.
It follows that all virtual $\g$-modules of the form $[U] - [U^{\dag}] \in R_{0}(\g)$ must lie in the kernel of $\Delta$.

On the other hand, suppose that $U_{1}$ and $U_{2}$ are two ungraded $\g$-modules such that $\Delta U_{1} \cong \Delta U_{2}$ as $\g$-supermodules.  If we forget the $\Z_{2}$-gradings, then Proposition \ref{composition} above gives us
$U_{1}\oplus U_{1}^{\dag} \cong U_{2}\oplus U_{2}^{\dag}$. Working
in the ungraded representation group $R_{0}(\g)$, we have
$$[U_{1}] - [U_{2}] = [U_{2}^{\dag}] - [U_{1}^{\dag}] =
  - \bigl( [U_{1}] - [U_{2}] \bigr)^{\dag}.$$
It then follows that the virtual $\g$-module $[U_{1}] - [U_{2}]$ lies
in the $-1$ eigen\-space of the involution $\dag :  R_{0}(\g) \to R_{0}(\g)$. Therefore the kernel of $\Delta$ is the anti-conjugate
representation group $\Rac(\g)$.

A similar argument works for the kernel of the map $R_{\Z_{2}}(\g) \to R_{0}(\g)$ which forgets the $\Z_{2}$-grading. If $V$ is a $\g$-supermodule, then $V$ is clearly isomorphic to its parity reversal $V^{\Pi}$ when we ignore their gradings. Thus the virtual $\g$-supermodule $[V] - [V^{\Pi}]\in R_{\Z_{2}}(\g)$ maps to $0\in R_{0}(\g)$ when we forget its grading.
If $V_{1}$ and $V_{2}$ are two $\g$-supermodules which are isomorphic as ungraded $\g$-modules, then we must have $\Delta V_{1} \cong \Delta V_{2}$ as $\g$-supermodules, and by Proposition \ref{composition} we obtain
an isomorphism
$V_{1}\oplus V_{1}^{\Pi} \cong V_{2}\oplus V_{2}^{\Pi}$. Working with virtual $\g$-supermodules in $R_{\Z_{2}}(\g)$, we obtain
$$[V_{1}] - [V_{2}] = [V_{2}^{\Pi}] - [V_{1}^{\Pi}] =
  - \bigl( [V_{1}] - [V_{2}] \bigr)^{\Pi}.$$
Thus $[V_{1}]-[V_{2}]$ is in the $-1$ eigenspace of the involution $\Pi: R_{\Z_{2}}(\g)\to R_{\Z_{2}}(\g)$, and so the kernel of the forgetful map is the anti-dual representation ideal $\Rm(\g)$.
\end{proof}

Now that we know the kernels and images of the diagonal map $\Delta$ and the forgetful map $f$, we can knit them together into a single periodic exact sequence.

\begin{theorem}\label{exact-sequence}
	If $\g$ is a Lie superalgebra, then its representation rings, ideals, and groups form a periodic exact sequence:
	\begin{equation}\label{eq:exact1}
	\begin{diagram}
		R_{0}(\g) & \rTo^{\Delta} & \Rp(\g) & \rTo^{\pi} & SR(\g) \\
		\uTo<{\delta'} & & & & \dTo>{\delta} \\
		SR'(\g) & \lTo^{\pi'} & \Rsc(\g) & \lTo^{f} & 
		R_{\Z_{2}}(\g)
	\end{diagram}
	\end{equation}
	where the vertical maps are the compositions
	\begin{equation}\label{eq:delta}
	\delta: SR(\g) = R_{\Z_{2}}(\g)/\Delta R_{0}(\g) \twoheadrightarrow R_{\Z_{2}}(\g)/\Rp(\g)
	\xrightarrow{1-\Pi} \Rm(\g) \hookrightarrow R_{\Z_{2}}(\g)
	\end{equation}
	and
	\begin{equation}\label{eq:delta1}
	\delta' : SR'(\g) = R_{0}(\g) / f\,R_{\Z_{2}}(\g) \twoheadrightarrow R_{0}(\g)/\Rsc(\g)
	\xrightarrow{1-\dag} \Rac(\g) \hookrightarrow R_{0}(\g),
	\end{equation}
	and the maps $\pi$ and $\pi'$ are given by the diagrams
\begin{equation}\label{eq:commute1} 
\begin{diagram}
	\Rp(\g) & \rInto & R_{\Z_{2}}(\g) \\
	\dOnto>{/\mathrm{Im}\,\Delta} & \rdTo^{\pi} & \dOnto>{/\mathrm{Im}\,\Delta} \\ 
	\Rp(\g)/\Delta R_{0}(\g) & \rInto &SR(\g) 
\end{diagram} \qquad\qquad
\begin{diagram}
	\Rsc(\g) & \rInto & R_{0}(\g) \\
	\dOnto>{/\mathrm{Im}\,f} & \rdTo^{\pi'} & \dOnto>{/\mathrm{Im}\,f} \\ 
	\Rsc(\g)/f\, R_{\Z_{2}}(\g) & \rInto & SR'(\g) 
\end{diagram}
\end{equation}
	(The maps here are all additive group homomorphisms, not ring homomorphisms.)
\end{theorem}

\begin{proof}
First we recall Proposition \ref{image}, which tells us that $\Delta R_{0}(\g) \subset \Rp(\g)$ and $f\,R_{\Z_{2}}(\g) \subset \Rsc(\g)$. This justifies our use of the self-dual representation ideal and the self-conjugate representation group as the ranges of the maps $\Delta$ and $f$ in (\ref{eq:exact1}).
Second, we see immediately from the commutative diagrams (\ref{eq:commute1}) that
$\mathrm{Ker}\,\pi = \mathrm{Im}\,\Delta$
and $\mathrm{Ker}\,\pi' = \mathrm{Im}\,f$,
which establishes the exactness of (\ref{eq:exact1}) at the top center and bottom center nodes.

Next we move on to the top right and bottom left nodes of (\ref{eq:exact1}). By Lemmas \ref{anti-sym} and \ref{anti-sym-dag}, the anti-symmetrizing operators $1-\Pi$ and $1-\dag$ appearing in  (\ref{eq:delta}) and  (\ref{eq:delta1}) are isomorphisms.  So by the definitions (\ref{eq:delta}) and (\ref{eq:delta1}) of the connecting maps $\delta$ and $\delta'$, we see that
$$\mathrm{Ker}\,\delta = \mathrm{Ker}\,\bigl( R_{\Z_{2}}(\g)/\Delta R_{0}(\g) \to R_{\Z_{2}}(\g)/\Rp(\g)\bigr),$$
and
$$\mathrm{Ker}\,\delta' = \mathrm{Ker}\,\bigl( R_{0}(\g)/f\, R_{\Z_{2}}(\g) \to R_{0}(\g)/\Rsc(\g)\bigr).$$
However, the sequences of inclusions
$$\Delta R_{0}(\g) \subset \Rp(\g) \subset R_{\Z_{2}}(\g)
\quad\text{and}\quad
 f\, R_{\Z_{2}}(\g) \subset \Rsc(\g) \subset R_{0}(\g)
$$
induce short exact sequences of the coset spaces
$$0\to\Rp(\g) / \Delta R_{0}(\g) \to
  R_{\Z_{2}}(\g) / \Delta R_{0}(\g) \to
  R_{\Z_{2}}(\g) / \Rp(\g) \to 0
$$
and
$$
 0\to\Rsc(\g) / f\,R_{\Z_{2}}(\g) \to
  R_{0}(\g) / f\, R_{\Z_{2}}(\g) \to
  R_{0}(\g) / \Rsc(\g)\to 0,
$$
respectively. So, we see that
$$\mathrm{Ker}\,\delta = 
\Rp(\g) / \Delta R_{0}(\g) = \mathrm{Im}\,\pi,\quad\mathrm{Ker}\,\delta' = \Rsc(\g) / f\,R_{\Z_{2}}(\g) =
\mathrm{Im}\,\pi'.$$

Finally, we consider the bottom right and top left nodes
of (\ref{eq:exact1}). From (\ref{eq:delta}) and (\ref{eq:delta1}),
we see that $\mathrm{Im}\,\delta = \Rm(\g)$ and $\mathrm{Im}\,\delta' = \Rac(\g)$. However, by Proposition \ref{kernel}, these are precisely the kernels of $f$ and $\Delta$,
respectively.  The periodic sequence (\ref{eq:exact1}) is therefore exact.
\end{proof}

\subsection{Examples}

\begin{example}
If $\g = \g_{0}$ is a Lie algebra, then the ungraded representation group is the usual non-supersymmetric representation ring $R_{0}(\g_{0}) = R(\g_{0})$. A $\g_{0}$-supermodule is simply a pair of $\g_{0}$-modules, and so the graded representation ring is two copies $R_{\Z_{2}}(\g_{0}) = R(\g_{0}) \oplus R(\g_{0})$. The parity reversal operator interchanges the two copies, and
thus a self-dual $\g_{0}$-supermodule is just two copies of the same ungraded $\g_{0}$-module. In other words, the self-dual representation ideal is the diagonal $\Rp(\g_{0}) = \Delta R_{0}(\g_{0}) \subset R_{\Z_{2}}(\g)$, and thus the diagonal map 
$\Delta: R_{0}(\g_{0}) \to \Rp(\g_{0})$ is an isomorphism.  The super representation ring is then the quotient $SR(\g_{0}) \cong R_{\Z_{2}}(\g_{0} ) / \Rp(\g_{0}) \cong R(\g_{0})$, given by formal differences of pairs of $\g_{0}$-modules. Since $\g_{0}$ has no odd component, every ungraded $\g_{0}$-module
is automatically self conjugate, so we have $\Rsc(\g_{0}) = R(\g_{0})$. Finally,
the forgetful map $f: R_{\Z_{2}}(\g_{0}) \to R_{0}(\g_{0})$, which adds the odd and even components, is onto, and thus the ungraded super representation group $SR'(\g_{0})$ vanishes.  All together, the periodic exact sequence (\ref{eq:exact1}) becomes:
$$
\begin{diagram}[tight,width=4em]
R(\g_{0}) & \rTo^{\cong}_{\Delta} & R(\g_{0}) & \rTo^{0}_{\pi} & R(\g_{0}) & \\
\uTo>{\delta'} & & & & \dInto<{\delta}>{1 \oplus -1} \\
0 & \lTo^{\pi'} & R(\g_{0}) & \lOnto^{f}_{+} & R(\g_{0}) \oplus R(\g_{0})
\end{diagram}
$$
\end{example}

\begin{example}
Consider the strange Lie superalgebra $\g = \mathfrak{q}(1)$ (see \cite{BL}), which has one even generator $H$ and one odd generator $Q$, with brackets
$$[H,H] = [H,Q] = 0, \qquad [Q,Q] = 2H.$$
In other words, we have $Q^{2} = H$ and all other brackets vanish. This is the simplest nontrivial Lie superalgebra, and it arises naturally in many supersymmetric settings. For instance, in \cite{Ko} and \cite{La}, both Kostant and the author implicitly use the superalgebra $\mathfrak{q}(1)$, where $Q = \dirac$ is Kostant's cubic Dirac operator and $H = \dirac^{2}$ is its square. A more familiar example is the de Rham complex, with $Q = d + d^{*}$ and $H = dd^{*} + d^{*}d$ the Laplacian.

We now consider complex representations of $\mathfrak{q}(1)$.
The irreducible $\mathfrak{q}(1)$-supermodules are the trivial representations $\mathrm{I} = \C \oplus 0$ and $\Pi = 0 \oplus \C$ (so named because the tensor product with $\mathrm{I}$ and $\Pi$ are the identity and parity reversal maps, respectively), as well as $L_{\lambda} = \C \oplus \Pi\C$ for each complex $\lambda \neq 0$, with action given by
$$H \mapsto \begin{pmatrix} \lambda & 0 \\ 0 & \lambda \end{pmatrix}, \quad
  Q \mapsto \begin{pmatrix} 0 & \sqrt{\lambda} \\ \sqrt{\lambda} & 0 \end{pmatrix}.$$
Every nontrivial irreducible $\mathfrak{q}(1)$-supermodule is then isomorphic to $L_{\lambda}$ for some complex $\lambda \neq 0$. These $L_{\lambda}$ for $\lambda \neq 0$ are all clearly self-dual, and we obtain, additively, the free abelian groups
$$R_{\Z_{2}}\bigl(\mathfrak{q}(1) \bigr) \cong \Z[\C^{*},\mathrm{I},\Pi], 
\qquad R_{+}\bigl(\mathfrak{q}(1) \bigr) \cong \Z[\C],$$
where $R_{+}(\mathfrak{q}(1))$ is generated by $[L_{\lambda}]$ for all $\lambda\in C$,
including $L_{0} = \mathrm{I} \oplus \Pi$. (Note that we are discussing only the additive structure here. We will consider the ring structure when we revisit this example in Section \ref{canonical-bases}.)
The irreducible ungraded $\mathfrak{q}(1)$-modules are all one-dimensional, and we denote by $\C_{\mu}$ for each $\mu\in\C$ the representation on which $H$ and $Q$ act by multiplication by $\mu^{2}$ and $\mu$, respectively. Conjugation flips the sign of the action of $Q$, giving $\C_{\mu}^{\dag} = \C_{-\mu}$. We thus have
$$R_{0}\bigl(\mathfrak{q}(1) \bigr) \cong \Z[\C], \qquad
  \Rsc\bigl(\mathfrak{q}(1) \bigr) \cong \Z[\C],$$
where $R_{0}(\mathfrak{q}(1))$ is generated by the classes $[\C_{\mu}]$ for $\mu\in \C$, while $\Rsc(\mathfrak{q}(1))$ is generated by the classes $[\C_{0}]$
and $[\C_{\sqrt{\lambda}}] + [\C_{-\sqrt{\lambda}}]$ for $0\neq\lambda\in \C$.
In this case, the diagonal map and the forgetful map are given by
$$\Delta : \bigl[\C_{\mu}\bigr] \mapsto \bigl[L_{\mu^{2}}\bigr], \qquad
  f : \bigl[L_{\lambda}\bigr] \mapsto \bigl[\C_{\sqrt{\lambda}}\bigr] + \bigl[\C_{-\sqrt{\lambda}}\bigr]$$
for $\mu,\lambda\neq 0$, while for the trivial representations, the diagonal map is
$\Delta [\C_{0}] = [\mathrm{I}] + [\Pi] = [L_{0}]$ and the forgetful map is $f[\mathrm{I}] = f[\Pi] = [\C_{0}]$.
This gives us 
\begin{align*}
	SR\bigl( \mathfrak{q}(1) \bigr) &= \Z[\mathrm{I},\Pi] / (\Pi = -\mathrm{I}) \cong \Z, \\
  SR' \bigl( \mathfrak{q}(1) \bigr) &= \Z[ \C^{*} ] / \bigl( [-\mu] =-[\mu]\bigr)
  \cong \Z[\C^{*}],\end{align*}
where $SR'(\mathfrak{q}(1))$ is generated by a choice of $\sqrt{\lambda}$ for each $\lambda \neq 0$.
The periodic exact sequence (\ref{eq:exact1}) now becomes
$$\begin{diagram}[tight,width=4em]
\Z[\C] & \rOnto_{\Delta}^{[\mu] \mapsto [\mu^{2}]} & \Z[\C] & \rTo^{0}_{\pi} & \Z \\
\uInto>{\delta'}<{[\lambda] \mapsto [\sqrt{\lambda}] - [-\sqrt{\lambda}]} & & & & \dInto<{\delta}>{\mathrm{I} \mapsto \mathrm{I}-\Pi} \\
\Z[ \C^{*} ] & \lTo^{0}_{\pi'} & \Z[\C] & \lOnto^{\substack{\;[\lambda]\mapsfrom [\lambda]\\  [0] \mapsfrom \mathrm{I},\Pi}}_{f} & \Z[\C^{*},\mathrm{I},\Pi]
\end{diagram}$$
\end{example}

\subsection{Further structure}
To help the reader keep track of the many variants of the
representation ring used in the proof of Theorem \ref{exact-sequence}, as well as the maps connecting them, we give here a complete diagram showing all of the groups and homomorphisms that we have thus far considered.
$$
\begin{diagram}[hug,tight,width=3.5em]
	\Rsc(\g)/\mathrm{Im}\,f & \lOnto &\Rsc(\g) & & \Rp(\g) & \rOnto & \Rp(\g)/\mathrm{Im}\,\Delta\\
	\dInto &\ldTo<{\pi'}& \dInto &  \luTo<{f\quad\qquad} \ruTo<{\quad\qquad\Delta} & \dInto & \rdTo>{\pi} & \dInto \\
	SR'(\g)& \pile{\lOnto \\ \rTo_{\delta'} \\ } &R_{0}(\g) & \pile{\rTo<{\Delta}\\ \lTo>{f}} & R_{\Z_{2}}(\g) & \pile{\rOnto \\ \lTo_{\delta}} & SR(\g) \\
	\dOnto & \ldOnto &\uInto &  & \uInto & \rdOnto & \dOnto \\
	R_{0}(\g) / \Rsc(\g) & \rTo^{\cong}_{1-\dag}&\Rac(\g)& & \Rm(\g) & \lTo^{\cong}_{1-\Pi} & R_{\Z_{2}}(\g) / \Rp(\g) 
\end{diagram}
$$
By examining the vertical sequences on the left side and right side of this diagram, we obtain useful decompositions of $SR(\g)$ and $SR'(\g)$.

\begin{proposition}\label{negative}
Both the parity reversal operator $\Pi : R_{\Z_{2}}(\g) \to R_{\Z_{2}}(\g)$ and the conjugation operator
$\dag : R_{0}(\g) \to R_{0}(\g)$ descend to involutions of
the super representation ring $SR(\g)$ and the ungraded super representation group $SR'(\g)$
satisfying $[V^{\Pi}] = - [V]$ and $[U^{\dag}] = -[U]$, respectively.
Furthermore, we have short exact sequences
\begin{equation}\label{eq:decompose-SR}
0\to \Rp(\g)/\Delta R_{0}(\g) \to SR(\g) \to R_{\Z_{2}}(\g)/\Rp(\g) \to 0
\end{equation}
and
\begin{equation}\label{eq:decompose-SR1}
0 \to \Rsc(\g) / f\,R_{\Z_{2}}(\g) \to SR'(\g) \to R_{0}(\g)/\Rsc(\g) \to 0,
\end{equation}
where the subgroups on the left are entirely $2$-torsion, and
the quotients on the right are freely generated.
\end{proposition}

\begin{proof}
We recall from Proposition \ref{image} that the parity reversal operator $\Pi$ fixes the image of the diagonal map $\Delta R_{0}(\g)$ and that the conjugation operator $\dag$ fixes the image
of the forgetful map $f\,R_{\Z_{2}}(\g)$. It follows that $\Pi$ descends to the super representation ring $SR(\g) = R_{\Z_{2}}(\g) / \Delta R_{0}(\g)$ and that $\dag$ descends to $SR'(\g) = R_{0}(\g) / f\,R_{\Z_{2}}(\g)$.

If $V$ is a $\g$-supermodule, then we recall from Proposition \ref{composition} that $\Delta V \cong V \oplus V^{\Pi}$. Since
the image of the diagonal map vanishes in $SR(\g)$, we see that $[V\oplus V^{\Pi}] =0\in SR(\g)$. Likewise, Proposition \ref{composition} also tells us that if $U$ is an ungraded $\g$-module, then $\Delta U \cong U \oplus U^{\dag}$, and so
$[U\oplus U^{\dag}] = 0\in SR'(\g)$.

If $V$ is a $\g$-supermodule with $V\cong V^{\Pi}$, then in $SR(\g)$
we have $$2\,[V] = [V \oplus V] = [V \oplus V^{\Pi}] = 0.$$ Therefore,
the image of the self-dual representation ideal $\Rp(\g)$ in $SR(\g)$ must be entirely $2$-torsion. Similarly, if $U$ is an ungraded $\g$-module with $U\cong U^{\dag}$, then in $SR'(\g)$
we have $2[U] = [U \oplus U] = [U \oplus U^{\dag}] = 0$, so the
image of the self-conjugate representation group $\Rsc(\g)$ in
$SR'(\g)$ must also be entirely $2$-torsion. In these cases,
we are forced to introduce the $2$-torsion to reconcile the
equations $[V^{\Pi}]=[V]$ or $[U^{\dag}]=[U]$ for a self-dual or self-conjugate representation with the identities $[V^{\Pi}]=-[V]$
and $[U^{\dag}]=-[U]$ we established above.

As for the quotients, the group $R_{\Z_{2}}(\g) / \Rp(\g)$ is generated additively by the classes of irreducible $\g$-supermodules $V$ such that
$V\ncong V^{\Pi}$, subject to the relations $[V^{\Pi}] = -[V]$. The resulting
group is then freely generated, with one generator chosen from each such
pair $[V]$, $[V^{\Pi}]$.  In the ungraded case, the group $R_{0}(\g) / \Rsc(\g)$ is generated additively by the classes of irreducible $\g$-modules $U$ such that
$U\ncong U^{\dag}$, subject to the relations $[U^{\dag}] = -[U]$. The resulting
group is then freely generated, with one generator chosen from each such
pair $[U]$, $[U^{\dag}]$.
\end{proof}

When the field $k$ is algebraically closed, we can apply Schur's Lemma
given in Section \ref{subsection-representation-ring}. In this case, an irreducible $\g$-supermodule $V$ which is isomorphic to its own parity reversal $V^{\Pi}$ must admit
a $\g$-equivariant parity reversing involution $\alpha : V \to V^{\Pi}$.  In light of Proposition \ref{prop-4}, every such self-dual $\g$-supermodule is therefore in the image of the diagonal map $\Delta$, and thus $\Delta R_{0}(\g) = R_{+}(\g)$.  It follows from the above Proposition that the super represention ring $SR(\g)$ has no 2-torsion
and is in fact given by $SR(\g) \cong R_{\Z_{2}}(\g) / R_{+}(\g)$. In terms of type M and type Q irreducibles (see Section \ref{parity-reversal}), we find that $SR(\g)$ is freely generated by the classes of irreducible supermodules $V$ of type M (i.e., $V\ncong V^{\Pi}$) with the identification $[V^{\Pi}] = -[V]$.
Using the methods of the following sections, we can also show that $f\,R_{\Z_{2}}(\g) = \Rsc(\g)$, and thus the ungraded super representation group $SR'(\g)$ likewise has no 2-torsion and is given by $SR'(\g) \cong R_{0}(\g) / \Rsc(\g)$. As a consequence of this discussion, when the field $k$ is not algebraically closed, the 2-torsion appearing in the extensions (\ref{eq:decompose-SR}) and (\ref{eq:decompose-SR1}) gives quantitative measures of the failure of Schur's Lemma.

\section{Clifford superalgebras and degree shifts}
\label{degree-shift}

In this and the following section, we specifically work over the complex numbers $k = \C$. We use Clifford algebras to construct degree-shifted
versions of all the representation rings, ideals, and groups defined in the previous sections, based on the Clifford algebraic construction of degree-shifted $K$-theory in \cite{ABS} and \cite{Ka}.

\subsection{Associative superalgebras}
In order to add Clifford algebras to the mix, we must first modify our discussion from the previous sections to include not only representation rings, ideals, and groups of Lie superalgebras, but also representation groups constructed from modules over associative (unital) superalgebras.  An associative superalgebra is simply a $\Z_{2}$-graded algebra, its unit being an even element. Given an associative superalgebra $A$, we can give $A$ the structure of a Lie superalgebra by defining the bracket
$$[a,b] = ab - (-1)^{|a|\,|b|} ba$$
for homogeneous elements $a$, $b$ with $\Z_{2}$-degrees $|a|$ and $|b|$, respectively. On the other hand, given a Lie superalgebra $\g$, we can construct its universal enveloping algebra
$$U(\g) := T^{\bullet}(\g) / \bigl(\,X Y - (-1)^{|X|\,|Y|}YX = [X,Y]\, \bigr),$$
where the tensor algebra $T^{\bullet}(\g) = \bigoplus_{k=0}^{\infty}T^{k}(\g)$ is the free associative algebra
generated by $\g$, and the relations are defined for homogeneous elements $X,Y\in \g$
with $\Z_{2}$-degrees $|X|$ and $|Y|$, respectively. The universal enveloping algebra is the unique associative superalgebra satisfying the property that every Lie superalgebra homomorphism $\g \to A$, where
$A$ is an associative superalgebra, factorizes uniquely through $U(\g)$. In particular, every representation $\g \to \End(V)$ lifts to an algebra homomorphism $U(\g) \to \End(V)$, and every algebra homomorphism $U(\g)\to\End(V)$ restricts to a representation $\g \to \End(V)$. As a consequence, the various representation rings, ideals, and groups defined for Lie superalgebras in the previous section generalize immediately to representation groups for associative superalgebras.

\begin{definition}
The super tensor product $A\tensor B$ of two superalgebras $A$ and $B$ is the superalgebra whose underlying vector space is the tensor product $A\otimes B$, but with the multiplication
$$ (a_{1} \tensor b_{1}) (a_{2} \tensor b_{2} )
 := (-1)^{|b_{1}|\,|a_{2}|} (a_{1} a_{2}) \tensor (b_{1}b_{2}),$$
for homogeneous elements $a_{1},a_{2}\in A$ and $b_{1},b_{2}\in B$.
\end{definition}

If both $A$ and $B$ are unital, then they inject into their super tensor product $A\tensor B$ as $A\tensor 1$ and $1\tensor B$, respectively. The elements of $A$ and $B$ supercommute with each other, and this multiplication differs from the standard tensor product of algebras in that the odd elements of $A$ anti-commute with the odd elements of $B$.
Given two $A$-modules, we can construct their exterior tensor product over $A \tensor A$, but to obtain an interior tensor product as an $A$-module, we must pull back via a comultiplication map $\Delta : A \to A \tensor A$. The universal enveloping algebra $U(\g)$ is in fact a Hopf superalgebra (see \cite{MM}), and so has such a comultiplication, which we use to construct the tensor product of $\g$-supermodules (see Proposition \ref{SR-tensor} below).

\subsection{Degree shifts}

\begin{definition}
	The $n$-times degree-shifted representation groups of a Lie superalgebra $\g$ are
	constructed from representations of $\g$ admitting $n$ additional
	supersymmetries. More precisely, for a fixed $n\geq 0$,
	the degree-shifted complex representation groups are given by
	\begin{align*}
	R_{\bullet}^{-n}(\g) &:= R_{\bullet}\bigl(U(\g)\tensor\Cl(n)\bigr) \\
	SR^{-n}(\g) &:= SR \bigl( U(\g) \tensor \Cl(n) \bigr),
	\end{align*}
	where the $\bullet$ stands for all the possible subscripts:
	$\Z_{2}$, $+$, $-$, $0$, $\mathrm{sc}$, $\mathrm{ac}$ for the graded, self-dual, anti-dual, ungraded, self-conjugate, and anti-conjugate representation groups, respectively.
\end{definition}

Here we are using the language of \cite{FHT,FHT1}, referring to the odd actions of the Clifford generators as supersymmetries.
Note that for $n > 0$, we obtain additive groups, not rings.
However, the tensor product of a $\Cl(p)$-supermodule and a $\Cl(q)$-supermodule is a $\Cl(p+q)$-supermodule (see \cite{ABS},\cite{LM}, or Proposition \ref{SR-tensor} below). As a consequence, if we take all gradings $n \geq 0$ simultaneously,
then $R_{\Z_{2}}^{\bullet}(\g)$
and $SR^{\bullet}(\g)$ become $\Z$-graded rings. In fact, due to the twofold
periodicity of the complex Clifford algebras, we actually obtain $\Z_{2}$-graded
rings, consisting of the even degrees $R_{\Z_{2}}^{0}(\g)$, $SR^{0}(\g)$ and the odd degrees $R_{\Z_{2}}^{1}(\g), SR^{1}(\g)$. We now establish these statements more precisely:

\begin{lemma}\label{periodicity}
	The degree-shifted complex representation groups of a Lie superalgebra $\g$ satisfy
	periodicity,
	\begin{equation*}
		R^{-n-2}_{\bullet}(\g) \cong R^{-n}_{\bullet}(\g),  \qquad
		SR^{-n-2}(\g) \cong SR^{-n}(\g),
	\end{equation*}
	as a consequence of the periodicity $\Cl(n+2) \cong \Cl(n) \otimes \C(2)$
	of complex Clifford algebras.
\end{lemma}

\begin{proof}
	Let $A = A_{0}\oplus A_{1}$ be an associative superalgebra.
	Let $e_{1}, e_{2}$ be two generators of $\Cl(2)$, and consider the element
	$\omega = i\,e_{1}\cdot e_{2} \in \Cl(2)$. We then have
	$$\omega^{2} = 1 \quad \text{ and } \quad \omega \cdot e_{i} = - e_{i}\cdot \omega \text{ for $i=1,2$.}
	$$
	Consider the map $f : A \tensor \Cl(2) \to A \otimes \Cl(2)$
	which is specified by
	\begin{equation}\label{eq:tensor-otimes}
		f( a_{0} \tensor 1 ) = a_{0} \otimes 1, \qquad
		f( a_{1} \tensor 1 ) = a_{1} \otimes \omega, \qquad
		f( 1 \tensor e_{i} ) = 1 \otimes e_{i},
	\end{equation}
	for odd and even elements $a_{0}\in A_{0}$ and $a_{1}\in A_{1}$ and
	for Clifford generators $e_{1}, e_{2}$.
	Since $\omega^{2} = 1$, the map $f$ preserves the multiplication on $A$,
	and it clearly preserves the multiplication on $\Cl(2)$.
	Furthermore, since $\omega$ anti-commutes with the Clifford generators $e_{i}$,
	we have
	\begin{equation*}
		f(a_{1} \tensor 1)\,f(1 \tensor e_{i}) + f(1\tensor e_{i})\,f(a_{1}\tensor 1 )
		= (a_{1} \otimes \omega \cdot e_{i})+ (a_{1} \otimes e_{i} \cdot \omega) = 0,
	\end{equation*}
	and so $f$ is a superalgebra homomorphism. If we exchange the super and
	standard tensor products $\tensor$
	and $\otimes$ in (\ref{eq:tensor-otimes}),
	then we obtain an inverse map
	$f^{-1} : A \tensor \Cl(2) \to A \otimes \Cl(2)$,
	which is likewise a superalgebra homomorphism.
	We therefore have an isomorphism $A \tensor \Cl(2) \cong A \otimes \Cl(2)$
	of superalgebras. Now, recalling the Clifford algebra
	isomorphism $\Cl(2) \cong \C(2)$ (see \cite{ABS}, \cite{LM}), we find that
	$A\tensor\Cl(2) \cong A \otimes \C(2)$. Therefore, by Morita equivalence,
	$A\tensor\Cl(2)$-modules correspond
	precisely to $A$-modules.
	
	So far we have established this result for the ungraded representation ring.
	However, we
	can regard $\Z_{2}$-graded modules over $A$ as ungraded modules over the
	superalgebra $A'$ generated by $A$ and the grading operator $(-1)^{F}$,
	so this result also holds for supermodules. The parity reversal and conjugation
	involutions commute with this construction, and so we further obtain isomorphisms
	of all the remaining representation rings.	
\end{proof}

\begin{proposition}\label{SR-tensor}
The tensor product of super vector spaces induces a product
$$R_{\Z_{2}}^{-n}(\g) \otimes R_{\Z_{2}}^{-m}(\g) \to R_{\Z_{2}}^{-n-m}(\g)$$
on the degree-shifted graded representation groups, which descends to a 
product 
$$SR^{-n}(\g) \otimes SR^{-m}(\g) \to SR^{-n-m}(\g)$$
on the degree-shifted super representation groups. Taking all degrees together,
$R_{\Z_{2}}^{\bullet}(\g)$ and $SR^{\bullet}(\g)$ are $\Z_{2}$-graded rings (and the degree homogeneous components are modules) over $R_{\Z_{2}}(\g)$ and $SR(\g)$, respectively. Furthermore, if $V$ and $W$ are $\g$-Clifford supermodules corresponding to classes in $R_{\Z_{2}}^{-n}(\g)$ and $R_{\Z_{2}}^{-m}(\g)$, respectively, then
$$V \otimes W \cong \begin{cases}
W \otimes V & \text{if at least one of $n,m$ is even,} \\
(W \otimes V)^{\Pi} & \text{if both $n,m$ are odd.}
\end{cases}
$$
Consequently, $[V][W] = (-1)^{nm}[W][V]$ in $SR^{\bullet}(\g)$,
and thus $SR^{\bullet}(\g)$ is supercommutative.
\end{proposition}

\begin{proof}
Let $V$ and $W$ be $\g$-supermodules carrying auxiliary actions of $\Cl(n)$ and $\Cl(m)$, respectively, which supercommute with the $\g$-action. The action on the tensor product $V \otimes W$ is then
\begin{alignat*}{2}
	X_{0} & \mapsto X_{0} \otimes 1 + 1 \otimes X_{0} & \qquad & \text{for even elements $X_{0}\in\g_{0}$}, \\
	X_{1} & \mapsto X_{1} \otimes 1 + (-1)^{F} \otimes X_{1} & \qquad & \text{for odd elements $X_{1}\in\g_{1}$}, \\
	e_{i} & \mapsto e_{i} \otimes 1 & \qquad & \text{for Clifford generators $e_{i} \in \Cl(n)$}, \\
	f_{j} & \mapsto (-1)^{F} \otimes f_{j} & \qquad & \text{for Clifford generators $f_{j}\in\Cl(m)$},
\end{alignat*}
where the $e_{i}$ and $f_{j}$ together generate the Clifford algebra $\Cl(n+m)$.
In other words, if $V$ is a $U(\g)\tensor\Cl(n)$-super\-module and $W$ is a $U(\g)\tensor\Cl(m)$-super\-module, then their tensor product $V\otimes W$ is a super\-module
over the super tensor product
$$\bigl( U(\g)\tensor\Cl(n) \bigr) \tensor  \bigl( U(\g) \tensor\Cl(m) \bigr)
\cong \bigl( U(\g) \tensor U(\g) \bigr) \tensor \Cl(n+m).$$
Such a representation then restricts to a $U(\g) \tensor \Cl(n+m)$-supermodule by taking the pullback with respect to the comultiplication map $\Delta: U(\g) \mapsto U(\g) \tensor U(\g)$.
We next show that this tensor product descends to the quotients
$$SR^{-n}(\g) \cong R_{\Z_{2}}^{-n}(\g) / i^{*} R_{\Z_{2}}^{-n-1}(\g)$$
which by Corollary \ref{karoubi-k} below give the degree-shifted super representation groups.
Namely, if $V$ or $W$ actually admits a $\Cl(n+1)$- or $\Cl(m+1)$-action, then their tensor product $V\otimes W$ clearly carries a $\Cl(n+m+1)$-action and therefore vanishes in the degree-shifted super representation group $SR^{-n-m}(\g)$. In both cases, we see that $R_{Z_{2}}^{\bullet}(\g)$ and $SR^{\bullet}(\g)$ are
$\Z$-graded rings over their degree zero components $R_{\Z_{2}}^{0}(\g) = R_{\Z_{2}}(\g)$ and $SR^{0}(\g) = SR(\g)$, respectively. The fact that they are
actually $\Z_{2}$-graded follows from the twofold periodicity of Lemma \ref{periodicity}.

If we reverse the order of the tensor product, $W \otimes V$ is once again a $U(\g) \tensor\Cl(n+m)$-supermodule.  It is isomorphic to $V \otimes W$ as a $\g$-supermodule, but the Clifford action may be different.
In particular, $V \otimes W$ is a $\Cl(n+m)$-supermodule with generators ordered $e_{1},\ldots e_{n},f_{1},\ldots,f_{m}$, while $W \otimes V$ is a $\Cl(n+m)$-supermodule with generators ordered $f_{1},\ldots,f_{m},e_{1},\ldots e_{n}$.
We observe that interchanging two Clifford generators of a $\g$-Clifford supermodule reverses its parity as follows: if the two Clifford generators are denoted by $e$ and $f$, then the operator $\alpha = (-1)^{F}\,2^{-1/2}\,(e + f)$ is an odd
involution which satisfies $e \alpha = \alpha f$ and $f \alpha = \alpha e$, but which commutes with the actions of $\g$ and the other Clifford generators.
Here, the permutation of the Clifford generators required to go from $V \otimes W$ to $W \otimes V$ takes $nm$ swaps, and thus $V \otimes W \cong W \otimes V$ if $nm$ is even and $V \otimes W \cong (W \otimes V)^{\Pi}$ if $nm$ is odd.
Recalling from Proposition \ref{negative} that the parity reversal operator $\Pi$ acts by $-1$ on the super representation ring $SR^{\bullet}(\g)$, we obtain $[V][W] = (-1)^{nm}[W][V]$, and thus $SR^{\bullet}(\g)$ is a supercommutative superring.
\end{proof}

\subsection{Canonical bases}
\label{canonical-bases}
At this point it is instructive to give canonical bases for the degree-shifted graded, self-dual, and super representation rings in terms of the irreducible $\g$-supermodules.
Recall from Schur's Lemma that an irreducible supermodule either has no odd automorphisms, or else its odd automorphisms are all nonzero scalar multiples of an odd involution.  These two classes of irreducibles are called type M (for matrix) and type Q (for queer), respectively (see \cite{BL}).  The parity reversal $M^{\Pi}$ of an irreducible of type M is then a distinct irreducible of type M, and so the type M irreducibles come in pairs, but on the other hand we have $Q^{\Pi} \cong Q$.
In degree $0$, we have the following additive bases over the integers for the graded, self-dual, and super representation rings:
\begin{align*}
	R_{\Z_{2}}(\g) &= \bigl\langle [M],\, [Q] \bigr\rangle \\
	\Rp(\g) &= \bigl\langle [M] +[M^{\Pi}],\, [Q] \bigr\rangle \\
	SR(\g) &= \bigl\langle [M] \bigr\rangle\,/\,\bigl( [M^{\Pi}] = -[M] \bigr)
\end{align*}
for all irreducible $\g$-supermodules $M$, $Q$ of types M, Q, respectively.

In degree $1$, we must consider pairs $(V,e)$ consisting of a $\g$-super\-module $V$ with a Clifford action $e$. If $M$ is an irreducible $\g$-supermodule of type M, then $M \oplus M^{\Pi}$ gives an irreducible of degree 1, because it in fact possesses \emph{two} supercommuting Clifford actions:
$$e_{1} = \begin{pmatrix}\phantom{-}0 & \mathrm{Id} \\ -\mathrm{Id} & 0 \end{pmatrix}, \qquad
 e_{2} = \begin{pmatrix}0 & i\,\mathrm{Id} \\ i\,\mathrm{Id} & 0 \end{pmatrix}.$$
Any choice $e = \cos\theta\,e_{1} + \sin\theta\,e_{2}$ for the Clifford action gives an
isomorphic pair $( M \oplus M^{\Pi}, e)$.  On the other hand, if $Q$ is an irreducible $\g$-supermodule of type Q, then we have only two possible Clifford actions: $e = \pm (-1)^{F} \alpha$, where $\alpha$ is the odd involution guaranteed by Schur's Lemma. Furthermore, $Q_{+}:= (Q,+(-1)^{F}\alpha)$ and $Q_{-} := (Q,-(-1)^{F}\alpha)$ are not (even) isomorphic to each other,
but $\alpha$ does give an odd isomorphism between them (i.e., $\alpha : Q_{+} \xrightarrow{\cong} Q^{\Pi}_{-}$).  As a consequence, in degree 1 the roles of type M and type Q irreducibles are reversed: the supermodule $M \oplus M^{\Pi}$ is now of type Q, while the supermodules $Q_{+}$ and $Q_{-}$ are now of type M.  The integral bases for the degree one representations rings are therefore
\begin{align*}
	R_{\Z_{2}}^{1}(\g) &= \bigl\langle [M\oplus M^{\Pi}],\, [Q_{+}], \, [Q_{-}] \bigr\rangle \\
	\Rp^{1}(\g) &= \bigl\langle [M \oplus M^{\Pi}],\, [Q_{+}] + [Q_{-}] \bigr\rangle \\
	SR^{1}(\g) &= \bigl\langle [Q_{+}] , \, [Q_{-}] \bigr\rangle\,/\,\bigl( [Q_{-}] = -[Q_{+}] \bigr)
\end{align*}
for all irreducible $\g$-supermodules $M$, $Q$ of types M, Q, respectively.

Putting these results together, we see that the full super representation ring $SR^{\bullet}(\g)$ is additively generated in degree 0 by half the type M irreducibles (one from each $M$, $M^{\Pi}$ pair) and in degree 1 by all the type Q irreducibles (or more precisely, one from each $Q_{+}$, $Q_{-}$ pair).  In contrast, each degree of the self-dual representation ring is the same additively; both $\Rp^{0}(\g)$ and $\Rp^{1}(\g)$ are generated by classes corresponding to half the type M and all the type Q irreducibles.  The same is true of the alternative version of the representation ring $R_{\Z_{2}}(\g) / \Rm(\g)$ in which we identify $[V] = [V^{\Pi}]$ (see Section \ref{parity-reversal}).  As a consequence, the true $K$-theoretic analogue to this alternative version of the representation ring is the full $\Z_{2}$-graded super representation ring $SR^{\bullet}(\g)$, which contains the same additive data, but has a distinct ($\Z_{2}$-graded) ring structure.

\begin{example}Let us consider once again the Lie superalgebra $\mathfrak{q}(1)$ which we introduced in Section \ref{les}.  It has irreducible supermodules $\mathrm{I}$ and $\Pi$ of type M, and $L_{\lambda}$ for nonzero $\lambda\in\C$ of type Q. As we saw before, in degree 0 we have the additive isomorphisms
\begin{align*}
R_{\Z_{2}}\bigl(\mathfrak{q}(1)\bigr) &\cong \Z[\C^{*},\mathrm{I},\Pi],\\
  \Rp(\mathfrak{q}\bigl(1)\bigr) &\cong \Z[\C],\\
  SR(\mathfrak{q}\bigl(1)\bigr) &\cong \Z[\mathrm{I},\Pi] / (\Pi = -\mathrm{I}) \cong \Z.
\end{align*}
Following our discussion above, in degree 1 we then have
\begin{align*}
R_{\Z_{2}}^{1}\bigl(\mathfrak{q}(1)\bigr) &\cong \Z[\mathrm{I}+\Pi,\C^{*}_{+},\C^{*}_{-}]
  \cong \Z[\C], \\
  \Rp^{1}(\mathfrak{q}\bigl(1)\bigr) &\cong \Z[\C], \\
  SR^{1}(\mathfrak{q}\bigl(1)\bigr) &\cong \Z[\C^{*}].
\end{align*}
(Note that these correspond to the groups we computed for $R_{0}$, $\Rsc$, and $SR'$ in Section \ref{les}.  This is not a coincidence, as we will show in the following section that expressions in terms of ungraded modules can be rephrased in terms of degree-shifted supermodules.) Here, the degree 0 component of the super representation ring is additively generated by the trivial supermodule $\mathrm{I}$ of type M, while the degree 1 component is generated by the supermodules $L_{\lambda}$ of type Q. The full super representation ring is then $SR^{\bullet}(\g) \cong \Z[\C]$, with $[0]$ in degree 0 and $[\lambda]$ for $\lambda\in\C^{*}$ in degree 1.  The product is given on these generators by $[0] [\lambda] = [\lambda] [0] = [\lambda]$ for all $\lambda\in\C$, and
$[\lambda_{1}] [\lambda_{2}] = 0$ for $\lambda_{1},\lambda_{2}\in\C^{*}$.  In contrast, we can consider the alternative version of the representation ring, constructed by identifying supermodules with their parity reversals, which is likewise $\Z[\C]$ as an additive group.  However, the product structure there is that of the group ring: $[\lambda_{1}][\lambda_{2}] = [\lambda_{1}+\lambda_{2}]$ for all $\lambda_{1}, \lambda_{2}\in \C$.
\end{example}

\section{The periodic exact sequence revisited}
\label{redux}

In this section we continue working over the complex numbers $k=\C$. At this point, we could simply work with the bases for the various representation rings described at the end of the previous section. However, we continue to work in full generality in order that we may carry over these results and their proofs to the real case in the following section.

The $\Z_{2}$-gradings for supermodules over a superalgebra $A$
can also be expressed in terms of Clifford algebras. The grading operator $(-1)^{F}$ acting on such a representation squares to the identity and anti-com\-mutes with the odd component of $A$. The operator $\sqrt{-1}\,(-1)^{F}$ therefore generates a complex Clifford algebra $\Cl(1)$, and so $\g$-supermodules can be viewed as
ungraded modules over the super tensor product $A \tensor \Cl(1)$. In other words, we have
\begin{equation}\label{eq:complex-grading}
	R_{\Z_{2}}^{-n}(\g) \cong R_{0}^{-n-1}(\g).
\end{equation}
Combining this with the Morita equivalence of Lemma~\ref{periodicity}, we obtain
the following:

\begin{proposition}\label{ungraded-shift}
	The ungraded degree-shifted representation groups of a Lie superalgebra $\g$
	are isomorphic to
	the degree-shifted representation groups of $\g$-super\-modules admitting
	one additional supersymmetry:
	$$R_{0}^{-n}(\g) \cong R_{\Z_{2}}^{-n-1}(\g).$$
\end{proposition}

\begin{proof}
	By (\ref{eq:complex-grading}), we have $R_{\Z_{2}}^{-n-1}(\g) \cong R_{0}^{-n-2}(\g)$, and by Lemma~\ref{periodicity} we have $R_{0}^{-n-2}(\g) \cong R_{0}^{-n}(\g)$.
\end{proof}

In the statement of Proposition~\ref{ungraded-shift} we once again use the language of \cite{FHT,FHT1}, referring to the odd actions of the Clifford generators as supersymmetries.
The isomorphism of Proposition \ref{ungraded-shift} is pivotal, since it allows us to rephrase statements about ungraded $\g$-modules in the language of degree-shifted $\g$-super\-modules. In this section, we examine this isomorphism carefully and recast our discussion from Sections \ref{lie-superalgebras} through \ref{les} in terms of the degree-shifted representation groups.

Recall the diagonal map $\Delta : R_{0}(\g) \to R_{\Z_{2}}(\g)$ defined in Section \ref{representation-rings}, which takes an ungraded $\g$-module $U$ to the graded $\g$-supermodule $\Delta U = U \oplus U$, with action given by (\ref{eq:deltar}). Now, consider the action of the Clifford generator $e_{1}$ on $\Delta U$ given by
$$c(e_{1}) = \begin{pmatrix} 0 & -\mathrm{Id} \\ \mathrm{Id} & 0\end{pmatrix}.$$
We then have $c(e_{1})^{2} = -\mathrm{Id}$, and furthermore $c(e_{1})$ commutes with the actions of the even elements and anti-commutes with the actions of the odd elements of $\g$.  The actions of $\g$ and $\Cl(1)$ therefore combine to give an action of the superalgebra $U(\g) \tensor \Cl(1)$ on $\Delta U$. If we consider the degree-shifted version, incorporating additional Clifford generators, we observe that the diagonal map $\Delta : R_{0}^{-n}(\g) \to R_{\Z_{2}}^{-n}(\g)$ in fact lifts to a map $\widetilde{\Delta} : R_{0}^{-n}(\g) \to R_{\Z_{2}}^{-n-1}(\g)$ such that $\Delta$ is the composition
$$R_{0}^{-n}(\g) \xrightarrow{\widetilde{\Delta}} R_{\Z_{2}}^{-n-1}(\g) \xrightarrow{i^{*}}
R_{\Z_{2}}^{-n}(\g),$$
where the map $i^{*}$ is induced by the inclusion
$i :\Cl(n)\hookrightarrow\Cl(n+1)$
and simply forgets the action of the extra Clifford generator.

To show that $\widetilde{\Delta}$ is the isomorphism of Proposition \ref{ungraded-shift}, we construct its inverse
$p_{0} : R^{-1}_{\Z_{2}}(\g) \to R_{0}(\g)$, which essentially projects onto
the even component of a $\g$-supermodule. We already encountered this
map in the proof of Proposition \ref{prop-4} in Section \ref{representation-rings}, and let us recall its construction.
Let $r : \g \to \End(V)$
be a $\Z_{2}$-graded representation of $\g$ on a super vector space
$V = V_{0} \oplus V_{1}$. We will now give an ungraded representation
$p_{0} r : \g \to \End(V_{0})$ on the even component $p_{0} V = V_{0}$.  Let $\alpha = (-1)^{F} e_{1}$ be the product of the grading operator with $e_{1}$. 
Since both the operators $(-1)^{F}$ and $e_{1}$ anti-commute with the odd elements of $\g$, their product $\alpha$ commutes with the odd elements of $\g$, and we obtain a $\g$-equivariant parity reversing involution $\alpha: V \to V^{\Pi}$.  Recall that the ungraded $\g$-action on $p_{0}V = V_{0}$ is then given by
\begin{equation}\label{eq:p0'}
	(p_{0}r)(X_{0}) = \bigl[r(X_{0}) \bigr]_{V_{0}} \quad\text{and}\quad
	(p_{0}r)(X_{1}) = \bigl[\alpha \, r(X_{1}) \bigr]_{V_{0}}
\end{equation}
for $X_{0}\in \g_{0}$ and $X_{1}\in \g_{1}$. Furthermore,
we can also define an ungraded $\g$-action on $p_{1}V := V_{1}$ by
the same method,
\begin{equation}\label{eq:p1}
	(p_{1}r)(X_{0}) = \bigl[r(X_{0}) \bigr]_{V_{1}} \quad\text{and}\quad
	(p_{1}r)(X_{1}) = \bigl[\alpha \, r(X_{1}) \bigr]_{V_{1}},
\end{equation}
and we observe that $\alpha$ gives an ungraded $\g$-module isomorphism $\alpha : p_{0}V \xrightarrow{\cong} p_{1} V$. Incorporating the additional Clifford generators of the degree-shifted version,
we get a map $p_{0}: R_{Z_{2}}^{-n-1}(\g) \to R_{0}^{-n}(\g)$.

\begin{proposition}\label{commutes}
	The following diagram commutes:
	\begin{equation}\label{eq:deltalift}
\begin{diagram}[vtrianglewidth=1.73em]
	R_{0}^{-n}(\g) & & \pile{ \rTo^{\widetilde{\Delta}} \\ \lTo_{p_{0}} } & & R_{\Z_{2}}^{-n-1}(\g) \\
	& \rdTo>{\Delta}&& \ldTo<{i^{*}} & \\
	&& R_{\Z_{2}}^{-n}(\g) &&
\end{diagram}
	\end{equation}
	and the maps $\widetilde{\Delta}$ and $p_{0}$ are inverses of each other.
\end{proposition}

\begin{proof}
We have already established in the above discussion that $\Delta = i^{*} \circ \widetilde{\Delta}$. In order to show that the maps $\widetilde{\Delta}$ and $p_{0}$ are inverses of each other, we verify that their composition $p_{0}\circ\widetilde{\Delta}$
is the identity map. We begin with an ungraded $\g$-module $U$ with
action $r : \g \to \End(U)$. As a super vector space, we have $\widetilde{\Delta} U = U \oplus U$, and projecting onto its even component gives us back our original vector space $p_{0}(\widetilde{\Delta} U) = U$. From our construction of $\widetilde{\Delta} U$, we observe that the map $\alpha = (-1)^{F} e_{1}$ is simply
$$\alpha = \begin{pmatrix}0&\mathrm{Id}\\\mathrm{Id}& 0\end{pmatrix}.$$ Combining equations (\ref{eq:deltar}) and (\ref{eq:p0'}), we see that  $p_{0}(\widetilde{\Delta} r) = r$, and thus $p_{0}(\widetilde{\Delta}U) = U$ as an ungraded $\g$-module. We now recall from the proof of Proposition \ref{prop-4} that $\widetilde{\Delta}(p_{0}V) \cong V$ for any $\g$-supermodule $V$, and thus $\widetilde{\Delta}$ and $p_{0}$ are inverse isomorphisms.
Finally, we verify the composition
$$\Delta \circ p_{0} = i^{*} \circ \widetilde{\Delta} \circ p_{0} = i^{*}$$
since $\widetilde{\Delta} \circ p_{0}$ is the identity map. 
\end{proof}

\begin{corollary}\label{karoubi-k}
Each degree-shifted super representation group of a Lie superalgebra is isomorphic to a quotient
$$SR^{-n}(\g)\cong R^{-n}_{\Z_{2}}(\g)\,/\,i^{*} R^{-(n+1)}_{\Z_{2}}(\g),$$
where the map $i^{*} : R^{-(n+1)}_{\Z_{2}}(\g) \to R^{-n}_{\Z_{2}}(\g)$
forgets the action of the $(n+1)$-st Clifford generator.
\end{corollary}

\begin{proof}
	This result follows immediately from the definition (\ref{eq:SR}) of the
	super representation ring, the commutative diagram (\ref{eq:deltalift}),
	and the Clifford algebra identity $\Cl(n) \tensor \Cl(1) \cong \Cl(n+1)$.
\end{proof}

\begin{corollary}
If $\g = \g_{0}$ is a Lie algebra, then its representation ring
$R^{\bullet}(\g_{0})$ and super representation ring $SR^{\bullet}(\g_{0})$ agree in all degrees,
$$SR^{-n}(\g_{0}) \cong R^{-n}(\g_{0}),$$
and in particular for the trivial Lie algebra $\g_{0} = \C$ we have
$$SR^{-n}(\C) = SR(\Cl(n)) \cong K^{-n}(\mathrm{pt}) \cong \tilde{K}(S^{n}).$$
\end{corollary}

\begin{lemma}\label{dag-pi}
	The isomorphism of Proposition \ref{ungraded-shift} intertwines the
	conjugation operator and the parity reversal operator. In other words,
	the following diagram commutes:
	$$\begin{diagram}
	R_{0}^{-n}(\g) & \pile{ \rTo^{\widetilde{\Delta}} \\ \lTo_{p_{0}} } & R_{\Z_{2}}^{-n-1}(\g) \\
	\dTo<{\dag} & & \dTo>{\Pi} \\
	R_{0}^{-n}(\g) & \pile{ \rTo^{\widetilde{\Delta}} \\ \lTo_{p_{0}} } & R_{\Z_{2}}^{-n-1}(\g)
	\end{diagram}$$
\end{lemma}
\begin{proof}
Given a $U(\g) \otimes \Cl(n+1)$-supermodule $V= V_{0}\oplus V_{1}$, let $\alpha = (-1)^{F} e_{n+1}$ be the product of the grading operator and the action of the last Clifford generator. Then $\alpha : V \to V^{\Pi}$ is a parity reversing involution of $V$.  Applying the projection maps $p_{0}$ and $p_{1}$ given by
(\ref{eq:p0'}) and (\ref{eq:p1}), respectively, we recall that $p_{0}V \cong p_{1}V$, where the isomorphism between these $\g$-modules is given by $\alpha : V_{0} \to V_{1}$. Reversing the parity of $V$ interchanges $V_{0}$ and $V_{1}$. Equivalently, it flips the sign of the grading operator $(-1)^{F}$, which in turn flips the sign of $\alpha$.  We then see from (\ref{eq:p0'}) and (\ref{eq:p1}) that $p_{0}(V^{\Pi}) = (p_{1}V)^{\dag}$ and that $p_{1}(V^{\Pi}) = (p_{0}V)^{\dag}$, since we are both swapping $p_{0}$ with $p_{1}$ and flipping the sign on the action of the odd component $\g_{1}$ of $\g$. It follows that $p_{0}(V^{\Pi}) = (p_{1} V)^{\dag} \cong (p_{0}V)^{\dag}$.

For the other direction, we use the fact from Proposition \ref{commutes} that the maps $\widetilde{\Delta}$ and $p_{0}$ are inverses of each other to obtain the desired composition:
$$\widetilde{\Delta} \circ \dag = \widetilde{\Delta} \circ \dag \circ p_{0} \circ \widetilde{\Delta} = \widetilde{\Delta} \circ p_{0} \circ \Pi \circ \widetilde{\Delta}
= \Pi \circ \widetilde{\Delta}$$
as a consequence of the identity $\dag \circ p_{0} = p_{0} \circ \Pi$ shown above.
\end{proof}

We can now reprise the discussion of Section \ref{les}, phrased entirely in terms of $\g$-supermodules and degree shifts. In
particular, we use the restriction
map $i^{*}: R_{\Z_{2}}^{-n-1}(\g) \to R_{\Z_{2}}^{-n}(\g)$ in place of the
diagonal map  $\Delta: R_{0}^{-n}(\g) \to R_{\Z_{2}}^{-n}(\g)$ and the parity
reversal operator $\Pi : R_{\Z_{2}}^{-n-1}(\g) \to R_{\Z_{2}}^{-n-1}(\g)$ in place of the conjugation operator $\dag : R_{0}^{-n}(\g) \to R_{0}^{-n}(\g)$.
In particular, our periodic exact sequence (\ref{eq:exact1}) becomes a periodic long exact sequence of degree-shifted representation groups.

\begin{proposition}\label{kernel-image}
Let $i^{*} : R_{\Z_{2}}^{-n-1}(\g) \to R_{\Z_{2}}^{-n}(\g)$ be the restriction
map induced by the inclusion $\Cl(n) \hookrightarrow \Cl(n+1)$. Then
$\mathrm{Ker}\,i^{*} = \Rm^{-n-1}(\g)$ and $\mathrm{Im}\,i^{*} \subset \Rp^{-n}(\g).$
\end{proposition}

\begin{proof}
We can apply Proposition \ref{commutes} to translate statements about the restriction map $i^{*}$ into statements about the diagonal map $\Delta$, and then this result is an immediate consequence of Propositions \ref{image} and \ref{kernel}. In particular, we have
$$\mathrm{Im}\,i^{*} = \mathrm{Im}(\Delta\circ p_{0}) = \mathrm{Im}\,\Delta \subset \Rp^{-n}(\g)$$
and
$$\mathrm{Ker}\,i^{*} = \mathrm{Ker}(\Delta \circ p_{0}) = \widetilde{\Delta}(\mathrm{Ker}\,\Delta ) = \widetilde{\Delta}(\Rac^{-n}(\g))
= \Rm^{-n-1}(\g),$$
where we apply Lemma \ref{dag-pi} to obtain the final equality.
\end{proof}

\begin{theorem}\label{exact-sequence1}
	If $\g$ is a Lie superalgebra, then its degree-shifted complex representation
	groups form a periodic long exact sequence:
	\begin{equation}\label{eq:exact}\begin{CD}
		\Rp^{0}(\g) @>\pi>> SR^{0}(\g) @>\delta>> R_{\Z_{2}}^{0}(\g) \\
		@Ai^{*}AA          @.         @Vi^{*}VV \\
		R^{1}_{\Z_{2}}(\g) @<\delta<< SR^{1}(\g) @<\pi<< \Rp^{1}(\g)
	\end{CD}\end{equation}
	where $\delta$ and $\pi$ are given by (\ref{eq:delta}) and
	(\ref{eq:commute1}) from Theorem
	\ref{exact-sequence}, with the appropriate degree shifts.
\end{theorem}

\begin{proof}
The periodicity of this sequence comes from the 2-fold periodicity of the degree-shifted representation groups shown in Lemma \ref{periodicity}. The maps preserve the degree except for the connecting maps $i^{*}$, which eliminate one Clifford generator and therefore raises the degree by 1. From Proposition \ref{kernel-image}, the image of $i^{*}$ lives in $i^{*}R_{\Z_{2}}^{-n-1}(\g)\subset R_{+}^{-n}(\g)$, and since $\pi$ is the composition
$$\pi : \Rp^{-n}(\g) \twoheadrightarrow \Rp^{-n}(\g)/i^{*}R_{\Z_{2}}^{-n-1}(\g)
\hookrightarrow R_{\Z_{2}}^{-n}(\g) / i^{*}R_{\Z_{2}}^{-n-1}(\g)
= SR^{-n}(\g),$$
we see that $\mathrm{Ker}\,\pi = \mathrm{Im}\,i^{*}$, making the sequence (\ref{eq:exact}) exact at the top left and bottom right corners.
The top center and bottom center nodes correspond to the top right node of
the exact sequence (\ref{eq:exact1}), and we recall from our proof of
Theorem \ref{exact-sequence} that $\mathrm{Im}\,\pi = \mathrm{Ker}\,\Delta$.
Finally, for the top right node, shifting the result of
Proposition \ref{kernel-image} by one degree gives us
$\mathrm{Ker}\,i^{*} = \Rm^{-n}(\g) \subset R_{\Z_{2}}^{-n}(\g)$. Since
$\delta$ is the composition
\begin{equation}\begin{split}
	\delta: SR^{-n}(\g) &= R_{\Z_{2}}^{-n}(\g)/i^{*} R_{\Z_{2}}^{-n-1}(\g) \\
	&\twoheadrightarrow R_{\Z_{2}}^{-n}(\g)/\Rp^{-n}(\g) 
	\xrightarrow{1-\Pi} \Rm^{-n}(\g) \hookrightarrow R_{\Z_{2}}^{-n}(\g),
\end{split}\end{equation}
where $1-\Pi$ is an isomorphism, we have $\mathrm{Im}\,\delta = \Rm^{-n}(\g) = \mathrm{Ker}\,i^{*}$, so the sequence (\ref{eq:exact}) is exact.
\end{proof}

In fact, since $\C$ is algebraically closed, we can apply Schur's Lemma (see Section \ref{subsection-representation-ring}).  As a consequence, any degree-shifted complex $\g$-super\-module $V$ which is isomorphic to its own parity reversal $V^{\Pi}$ admits a $\g$-equivariant parity reversing involution $\alpha : V \to V^{\Pi}$.  It follows from Propositions \ref{prop-4} and \ref{commutes} that
$$R_{+}^{-n}(\g) = \Delta R_{0}^{-n}(\g) = i^{*} R_{\Z_{2}}^{-n-1}(\g).$$
The vertical maps $i^{*}$ in the periodic long exact sequence (\ref{eq:exact}) are therefore surjections, and thus (\ref{eq:exact}) splits into two short exact sequences:
\begin{align*}
	0 \to SR^{0}(\g) \xrightarrow{\delta} R_{\Z_{2}}^{0}(\g)
	\xrightarrow{i^{*}}  R^{1}_{+}(\g) \to 0, \\
	0 \to SR^{1}(\g) \xrightarrow{\delta} R_{\Z_{2}}^{1}(\g)
	\xrightarrow{i^{*}}  R^{0}_{+}(\g) \to 0.
\end{align*}
From this, we see that the super representation ring and the degree-shifted self-dual representation group are essentially complementary inside of the graded representation ring.
As a consequence, if we consider the second sequence but substitute the
full graded representation ring $R_{\Z_{2}}^{0}(\g)$ for the self-dual representation ideal $R_{+}^{0}(\g)$, we can extend this exact sequence to
$$0 \to SR^{1}(\g) \to R^{-1}_{\Z_{2}}(\g) \xrightarrow{i^{*}} R^{0}_{\Z_{2}}(\g) \to SR^{0}(\g) \to 0$$
or equivalently, by Proposition \ref{commutes}, we have
$$0 \to SR^{1}(\g) \to R_{0}(\g) \xrightarrow{\Delta} R_{\Z_{2}}(\g) \to SR^{0}(\g) \to 0.$$
So, when working with the complex super representation ring, we can define
\begin{align*}
	SR^{0}(\g) &:= \mathrm{Coker}\,\Delta = \mathrm{Coker}\,i^{*}, \\
	SR^{1}(\g) &:= \mathrm{Ker}\,\Delta \cong \mathrm{Ker}\,i^{*}.
\end{align*}
Recalling Proposition \ref{kernel-image}, we find that as additive groups,
$$	SR^{0}(\g) \cong \Rm^{0}(\g), \qquad
	SR^{1}(\g) \cong \Rm^{1}(\g).$$
This simplification does not hold for the real super representation ring, as we show in the following section.

Finally, we demonstrate that the periodic long exact sequence of degree-shifted representation groups given by Theorem \ref{exact-sequence1} agrees with our periodic exact sequence of graded and ungraded representation groups
given by Theorem \ref{exact-sequence} from Section \ref{les}.

\begin{lemma}\label{restriction-equals-forgetful}
Under the isomorphism (\ref{eq:complex-grading}), the restriction map $i^{*}$ coincides with the forgetful map $f$, and the parity reversal operator $\Pi$ coincides with the conjugation operator $\dag$. In other words, the following diagrams commute:
$$
\begin{diagram}[vtrianglewidth=1.73em]
 && R_{\Z_{2}}^{-n-1}(\g) &&  \\
  & \ldTo<{i^{*}} & & \rdTo>{f} & \\
R_{\Z_{2}}^{-n}(\g) && \rTo^{\cong} && R_{0}^{-n-1}(\g)
\end{diagram} \qquad \qquad \qquad
\begin{diagram}
R_{\Z_{2}}^{-n}(\g) & \rTo^{\cong} & R_{0}^{-n-1}(\g) \\
\dTo<{\Pi} & & \dTo>{\dag} \\
R_{\Z_{2}}^{-n}(\g) & \rTo^{\cong} & R_{0}^{-n-1}(\g)
\end{diagram}
$$
and so the map (\ref{eq:complex-grading}) induces an isomorphism
$\Rsc^{-n-1}(\g) \cong \Rp^{-n}(\g).$
\end{lemma}
\begin{proof}
The key fact here is that the grading operator $(-1)^{F}$ can be treated as an additional Clifford generator. In particular, it commutes with all even actions and anti-commutes with all odd ones, and it satisfies $(-1)^{F} = \mathrm{Id}$.
Since we are working over the complex numbers, we then have $\bigl( \sqrt{-1}\,(-1)^{F} \bigr)^{2} = -\mathrm{Id}$, putting this Clifford generator in its usual form. Now, the restriction map $i^{*}$ forgets the action of the last Clifford generator, while the map $f$ forgets the $\Z_{2}$-grading. Each simply forgets one Clifford generator, and thus these maps are the same.

Recall from Remark~\ref{remark:conjugation} that $\g$-supermodules are automatically self-conjugate, but if we treat the grading operator $(-1)^{F}$ as an additional Clifford generator, then we must also
reverse the sign of the grading operator when we conjugate. However, reversing
the sign of the grading operator simply interchanges the even and odd components. Given a $\g$-supermodule $V\in R_{\Z_{2}}^{-n}(\g)$, let $V'\in R_{0}^{-n-1}(\g)$ be the corresponding ungraded $\g$-module under the isomorphism (\ref{eq:complex-grading}). We then have
$$(V')^{\dag} = \bigl( (V^{\dag})^{\Pi} \bigr)' \cong (V^{\Pi})',$$
and thus parity reversal indeed corresponds to conjugation.
\end{proof}

\begin{corollary}\label{equivalent}
	The ungraded and degree-shifted super representation groups,
	$$SR'(\g) := R_{0}(\g) / f\,R_{\Z_{2}}(\g) \qquad \text{and} \qquad
	SR^{1}(\g) := R_{\Z_{2}}^{-1}(\g) / i^{*}R_{\Z_{2}}^{0}(\g)$$
	are equivalent, and the periodic long exact sequences (\ref{eq:exact1}) and (\ref{eq:exact}) are isomorphic. 
\end{corollary}

Rephrasing our discussion following Theorem \ref{exact-sequence1} in terms of Corollary \ref{equivalent}, we find that in the complex case, the exact sequence (\ref{eq:exact1}) splits into two short exact sequences
\begin{align*}
	0 \to SR(\g) \xrightarrow{\delta} R_{\Z_{2}}(\g)
	\xrightarrow{f}  \Rsc(\g) \to 0, \\
	0 \to SR'(\g) \xrightarrow{\delta} R_{0}(\g)
	\xrightarrow{\Delta}  R_{+}(\g) \to 0.
\end{align*}
Extending the first of these, we then obtain the exact sequence
$$0 \to SR(\g) \to R_{\Z_{2}}(\g) \xrightarrow{i^{*}} R_{0}(\g) \to SR'(\g) \to 0$$
So, when working over the complex numbers, we can take
\begin{equation*}
	SR(\g)  := \mathrm{Ker}\,f, \qquad SR'(\g) := \mathrm{Coker}\,f
\end{equation*}
as our definition of the super representation ring and the ungraded super representation group.

\section{The real case}
\label{real}

The results of this paper also hold, with minor modifications, in the real case, where we consider real representations rather than complex ones and use real Clifford algebras $\mathrm{Cl}(p,q)$ in place of their complex counterparts for performing degree shifts. In this section, we describe the differences in the real case. Here, we denote all of the various real representation rings and groups by $R\R_{\bullet}(\g)$ and $\SRR(\g)$ in place of $R_{\bullet}(\g)$ and $SR(\g)$. The definitions of all these rings and groups, as well as the maps between them, work equally well for either real or complex $\g$-modules or supermodules; nowhere did we explicitly use $\sqrt{-1}$. For degree shifts,
we use the real Clifford algebra $\mathrm{Cl}(n)$ generated by the $n$ elements $\{e_{1}, \ldots, e_{n}\}$ subject to the Clifford relations (\ref{eq:clifford}), and we define
\begin{align*}
R\R_{\bullet}^{-n}(\g) &:= R\R_{\bullet}\bigl(U(\g) \tensor \mathrm{Cl}(n) \bigr), \\
\SRR^{-n} &:= \SRR\bigl(U(\g) \tensor \mathrm{Cl}(n) \bigr) \\
		  &\phantom{:}\cong R\R_{\Z_{2}}^{-n}(\g) / \Delta R\R_{0}^{-n}(\g) \\
		  &\phantom{:}\cong R\R_{\Z_{2}}^{-n}(\g) / i^{*} R\R_{\Z_{2}}^{-n-1}(\g).
\end{align*}
In fact, in the real case, we can consider the Clifford algebras $\mathrm{Cl}(p,q)$, with two sets of anti-commuting generators: $\{e_{1}, \ldots, e_{p}\}$ squaring to $-1$ and
$\{f_{1}, \ldots, f_{q}\}$ squaring to $+1$. In other words, we have the relations:
\begin{alignat*}{3}
e_{i}^{2} &= -1 \text{ for }i=1,\ldots,p,&\qquad e_{i} \cdot e_{j} &= - e_{j} \cdot e_{i} \text{ for $i\neq j$},\\
f_{k}^{2} &= +1\text{ for }k=1,\ldots,q, &\qquad f_{k} \cdot f_{l} &= - f_{l}\cdot f_{k} \text{ for $k\neq l$}, \\
&&e_{i }\cdot f_{k} &= - f_{k}\cdot e_{i} \text{ for all $i,k$}
\end{alignat*}
(i.e., $\mathrm{Cl}(p,q)$ is the Clifford algebra on $\R^{p,q}$ with the standard quadratic form of signature $(p,q)$). In particular, we have the Clifford algebra identifications $\mathrm{Cl}(p,q) \otimes \C \cong \Cl(p+q)$ and
$\mathrm{Cl}(n,0) = \mathrm{Cl}(n)$. Using these Clifford algebras, we obtain bidegree-shifted representation groups
\begin{align*}
	R\R_{\bullet}^{-p,-q}(\g) &:= R\R_{\bullet}\bigl(U(\g) \tensor \mathrm{Cl}(p,q)\bigr), \\
\SRR^{-p,-q} &:= \SRR\bigl(U(\g) \tensor \mathrm{Cl}(p,q) \bigr) \\
		  &\phantom{:}\cong R\R_{\Z_{2}}^{-p,-q}(\g) / \Delta R\R_{0}^{-p,-q}(\g) \\
		  &\phantom{:}\cong R\R_{\Z_{2}}^{-p,-q}(\g) / i^{*} R\R_{\Z_{2}}^{-p-1,-q}(\g),
\end{align*}
where we observe that the restriction map $i^{*}$ forgets the action of the last Clifford generator $e_{p+1}$ satisfying $e_{p+1}^{2} = -1$.

\subsection{Periodicities}
\label{real-periodicities}
In the real case, we have two distinct periodicities.  First, if we follow
the proof of Lemma \ref{periodicity} using real Clifford algebras, then in
order to obtain an element $\omega$ with $\omega^{2}=1$ that anti-commutes
with the two Clifford generators, we must take $\omega = e_{1} \cdot f_{1}\in\mathrm{Cl}(1,1)$. This then gives us the $(1,1)$-periodicity:
\begin{alignat*}{2}
	R\R^{-p-1,-q-1}_{\bullet}(\g) \cong R\R^{-p,-q}_{\bullet}(\g), \qquad
	SR_{\R}^{-p-1,-q-1}(\g) \cong SR_{\R}^{-p,-q}(\g),
\end{alignat*}
as a consequence of the real Clifford algebra periodicity
$$\mathrm{Cl}(p+1,q+1) \cong \mathrm{Cl}(p,q) \otimes \mathrm{Cl}(1,1)
\cong \mathrm{Cl}(p,q)\otimes \R(2).$$
Using this $(1,1)$-periodicity, we find that
\begin{equation}\label{eq:one-one}
	R\R^{p,q}_{\bullet}(\g) \cong R\R^{p-q}_{\bullet}(\g),  \qquad
	\SRR^{p,q}(\g) \cong \SRR^{p-q}(\g),
\end{equation}
and thus all of the bidegree-shifted real representation groups can be expressed in terms of the standard degree-shifted groups.

Second, we have the eightfold periodicity of the degree-shifted real representation groups:
\begin{equation*}
	R\R^{-n-8}_{\bullet}(\g) \cong R\R^{-n}_{\bullet}(\g), \qquad
	\SRR^{-n-8}(\g) \cong \SRR^{-n}(\g),
\end{equation*}
which follows from the eightfold periodicity
$$\mathrm{Cl}(n+8) \cong \mathrm{Cl}(n) \otimes \mathrm{Cl}(8)
\cong \mathrm{Cl}(n)\otimes \R(16)$$
of real Clifford algebras. For further information regarding the periodicities of real Clifford algebras, see \cite{A} or \cite{LM}.  The real version of Proposition \ref{SR-tensor} then tells us that the full graded representation ring $R\R_{\Z_{2}}^{\bullet}(\g)$ and the super representation ring $SR_{\R}^{\bullet}(\g)$ are both $\Z_{8}$-graded rings over their degree zero components $R\R_{\Z_{2}}(\g)$ and $SR_{\R}(\g)$, respectively, and that furthermore $SR_{\R}(\g)$ is supercommutatitive.

\subsection{The 24-term periodic exact sequence}
We now revisit Proposition \ref{ungraded-shift}, which was the keystone result of Section \ref{redux}.
Since we do not have the luxury of multiplying by $\sqrt{-1}$ when working over the real numbers, we note that the grading operator $(-1)^{F}$ on a super vector space squares to $+1$ and thus generates the real Clifford algebra $\mathrm{Cl}(0,1)$. As a consequence, the real version of equation (\ref{eq:complex-grading}) is the isomorphism
\begin{equation}\label{eq:real-grading}
	R\R_{\Z_{2}}^{-p,-q}(\g) \cong R\R_{0}^{-p,-q-1}(\g),
\end{equation}
which views the $\Z_{2}$-grading as simply one additional Clifford generator of type $(0,1)$.
Then, due to the $(1,1)$-periodicity of real Clifford algebras, we obtain
$$ R\R_{0}^{-p,-q}(\g) \cong R\R_{0}^{-p-1,-q-1}(\g)
   \cong R\R_{\Z_{2}}^{-p-1,-q}(\g),$$
or equivalently in terms of the degree rather than the bidegree, we have 
\begin{equation}\label{eq:real-grading-degree}
	R\R_{0}^{-n}(\g) \cong R\R_{\Z_{2}}^{-n-1}(\g),
\end{equation}
which agrees with the statement of Proposition \ref{ungraded-shift}. The rest of
the results of Section \ref{redux} up through Theorem \ref{exact-sequence1} continue to hold in the real case without modification, and we once again obtain a periodic long exact sequence, although now with an eightfold instead of twofold periodicity.

\begin{theorem}\label{24}
	If $\g$ is a Lie superalgebra, then its degree-shifted real representation
	groups form a periodic long exact sequence:
	\begin{equation}\label{eq:longexact}
		\cdots \to 
		R\R_{\Z_{2}}^{-n-1}(\g) \xrightarrow{i^{*}} R\R_{+}^{-n}(\g) \to
		\SRR^{-n}(\g)
		\to R\R_{\Z_{2}}^{-n}(\g) \xrightarrow{i^{*}} \cdots
	\end{equation}
	with periodicity of order $8$ in $n$ (i.e., the entire sequence repeats after 24 terms).
\end{theorem} 

\begin{example}
For the trivial Lie algebra $\g = \R$, we recover
$$\SRR^{-n}(\R) = \SRR(\mathrm{Cl}(n))
  \cong KO^{-n}(\mathrm{pt}) \cong \widetilde{KO}(S^{n}).$$
As in the complex case, the self-dual representation ideal is always  $R\R_{+}^{-n}(\R) = \Z$ for all $n$. The other representation rings in the
exact sequence take the values
\begin{align*}
	\SRR^{-n}(\R) &= \begin{cases}
	\Z & \text{for $n = 4k$,} \\
	\Z_{2} & \text{for $n = 8k + 1, 8k+ 2$,} \\
	0  & \text{otherwise}.
	\end{cases} \\
  R\R_{\Z_{2}}^{-n}(\R) &= \begin{cases}
    \Z \oplus \Z & \text{for $n = 4k$,} \\
    \Z & \text{otherwise}.
  \end{cases}
\end{align*}
The $\Z_{2}$-torsion terms appear since the restriction map $i^{*}$ in the exact sequence is multiplication by $2$ when $n = 8k + 1$ or $n = 8k+2$. All of these terms
then combine together to form the 24-term exact sequence given in Figure \ref{figure}.
\begin{figure}
\begin{center}
$$\begin{diagram}[hug]
		   & R\R_{+}^{-n}(\R) & & \SRR^{-n}(\R) & & R\R_{\Z_{2}}^{-n}(\R) \\
	n = 8: & \Z & \rTo_{0} & \Z & \rInto  & \Z\oplus\Z \\
	n = 7: & \Z & \rTo & 0 & \rTo \ldOnto(4,1) & \Z \\
 	n = 6: & \Z & \rTo & 0 & \rTo \ldTo(4,1)<{\cong\qquad} & \Z \\
	n = 5: & \Z & \rTo & 0 & \rTo \ldTo(4,1)<{\cong\qquad} & \Z \\
	n = 4: & \Z & \rTo_{0} & \Z & \rInto \ldTo(4,1)<{\cong\qquad} & \Z\oplus\Z \\
	n = 3: & \Z & \rTo & 0 & \rTo \ldOnto(4,1) & \Z \\
	n = 2: & \Z & \rOnto & \Z_{2} & \rTo^{0} \ldInto(4,1)<{\times 2\qquad}& \Z \\
	n = 1: & \Z & \rOnto & \Z_{2} & \rTo^{0} \ldInto(4,1)<{\times 2\qquad} & \Z \\
	n = 0: & \Z & \rTo_{0} & \Z & \rInto \ldTo(4,1)<{\cong\qquad} & \Z\oplus\Z \\
\end{diagram}$$
(The $n=0$ row corresponds to the $n=8$ row.)
\caption{The 24-term exact sequence for $\g=\R$ trivial.}
\label{figure}
\end{center}
\end{figure}
\end{example}

\subsection{The six-term periodic exact sequences}
On the other hand, the proof of Lemma \ref{restriction-equals-forgetful} works only in the complex case, and as a result, Corollary \ref{equivalent} no longer holds in the real case.  Surprisingly, we do still obtain a six-term periodic exact sequence for the real representation rings analogous to (\ref{eq:exact1}), but it is distinctly not isomorphic to our 24-term periodic long exact sequence (\ref{eq:longexact}). In order to build the real version of the periodic exact sequence (\ref{eq:exact1}), we must consider both the diagonal map $\Delta: R\R_{0}(\g) \to R\R_{\Z_{2}}(\g)$ and the forgetful map $f : R\R_{\Z_{2}}(\g) \to R\R_{0}(\g)$. If we express these maps in terms of the bidegree-shifted representation rings, then we recall that $\Delta$ is equivalent to the restriction map
\begin{equation}\label{eq:real-delta}
	i^{*} : R\R_{\Z_{2}}^{-p,-q}(\g) \to R\R_{\Z_{2}}^{-p+1,-q}(\g),
\end{equation}
while in light of (\ref{eq:real-grading}), the forgetful map is equivalent to
\begin{equation}\label{eq:real-f}
	f : R\R_{\Z_{2}}^{-p,-q}(\g) \to R\R_{\Z_{2}}^{-p,-q+1}(\g).
\end{equation}
In other words, the diagonal map $\Delta$, or equivalently the restriction map $i^{*}$, forgets one of the Clifford generator of type $(1,0)$, while the forgetful map $f$ forgets one of the Clifford generators of type $(0,1)$ (i.e., a $\Z_{2}$-grading).

These two maps (\ref{eq:real-delta}) and (\ref{eq:real-f}) are mirror images of each other if we interchange the type $(1,0)$ and type $(0,1)$ Clifford generators.  In general, given a Lie superalgebra $\g$, we can construct a new Lie superalgebra $\bar{\g}$ with the same underlying vector space, but with the new bracket
$$[X,Y]_{\bar{\g}} = (-1)^{|X|\,|Y|}[X,Y]_{\g},$$
for homogeneous elements $X,Y$ of $\Z_{2}$-degrees $|X|,|Y|$, respectively. The bracket on $\bar{\g}$ agrees with the bracket on $\g$ except for a change of sign when taking the bracket of two odd elements. Alternatively,
viewing a Lie superalgebra as consisting of a Lie algebra $\g_{0}$ paired with a $\g_{0}$-module $\g_{1}$ equipped with a $\g_{0}$-invariant symmetric bilinear form $\g_{1} \otimes \g_{1} \to \g_{0}$, taking the construction of $\bar{\g}$ simply flips the sign of the symmetric bilinear form. The same construction works for (unital) associative superalgebras, and we find that for real Clifford algebras, we have $\overline{\mathrm{Cl}(p,q)} = \mathrm{Cl}(q,p)$, giving us our mirror map.

\begin{proposition}\label{g-bar}
	The graded representation rings of the Lie superalgebras $\g$ and
	$\bar{\g}$ are isomorphic, but with opposite degrees:
	$$R\R_{\Z_{2}}^{-n}(\g) \cong R\R^{+n}_{\Z_{2}}(\bar{\g}).$$
\end{proposition}
\begin{proof}
For real Clifford algebras the isomorphisms
$$\mathrm{Cl}(p,q+2) \cong \mathrm{Cl}(q,p) \otimes \mathrm{Cl}(0,2)
\cong \mathrm{Cl}(q,p) \otimes \R(2)$$
(see \cite{LM}) induce the following isomorphisms of bidegree-shifted representation groups:
$$R\R_{\Z_{2}}^{-p,-q-1}(\g) 
\cong R\R_{0}^{-p,-q-2}(\g)
\cong R\R_{0}^{-q,-p}(\bar{\g})
\cong R\R_{\Z_{2}}^{-q,-p+1}(\bar{\g})$$
via an argument analogous to the proof of Lemma \ref{periodicity}. Applying (\ref{eq:one-one}) and putting $n=p-q-1$, the degree-shifted version of this isomorphism is then our desired result
$R\R_{\Z_{2}}^{n}(\g) \cong R\R_{\Z_{2}}^{-n}(\bar{\g}).$
\end{proof}

To describe the isomorphism of Proposition \ref{g-bar} explicitly, suppose that $V = V_{0}\oplus V_{1}$ is a $\g$-supermodule with action $r : \g \to \End(V)$.  We can define a new action $\bar{r} : \bar{\g} \to \End V$ by putting
\begin{equation*}
	\bar{r}(X_{0}) = r(X_{0}), \qquad
	\bar{r}(X_{1})|_{V_{0}} = + r(X_{1})|_{V_{0}}, \qquad
	\bar{r}(X_{1})|_{V_{1}} = - r(X_{1})|_{V_{1}}
\end{equation*}
for $X_{0}\in \g_{0}$ and $X_{1}\in \g_{1}$.  It is clear that $\bar{r}$ is indeed a representation of $\bar{\g}$, and since applying this
map twice gives us back our original $\g$-supermodule (and applying it twice
to a $\bar{\g}$-supermodule likewise acts by the identity), this map must be an isomorphism. We use this map to examine the real version of the definition (\ref{eq:SR1}) of the ungraded super representation group.

\begin{proposition}\label{restriction-forgetful-intertwine}
The isomorphism of Proposition \ref{g-bar} intertwines the restriction map $i^{*}$ and the forgetful map $f$ of equations (\ref{eq:real-delta}) and (\ref{eq:real-f}), respectively. In other words, the diagram
$$
\begin{diagram}
	R\R_{\Z_{2}}^{-n}(\g) & \rTo^{f} & R\R_{\Z_{2}}^{-n-1}(\g) \\
	\dTo<\cong & & \dTo>\cong \\
	R\R_{\Z_{2}}^{n}(\bar{\g}) & \rTo^{i^{*}} & R\R_{\Z_{2}}^{n+1}(\bar{\g})
\end{diagram}
$$
commutes, and therefore
$SR_{\R}'(\g) = R\R_{0}(\g) / f\,R\R_{\Z_{2}}(\g) \cong SR^{1}_{\R}(\bar{\g}).$
\end{proposition}

\begin{proof}
We leave to the reader the proof that the diagram commutes. As for the quotient,
we have the chain of isomorphisms:
\begin{equation*}\begin{split}
SR_{\R}'(\g) = R\R_{0}(\g) / f\,R\R_{\Z_{2}}(\g) &\cong R\R_{\Z_{2}}^{-1}(\g) / f\,R\R_{\Z_{2}}^{0}(\g) \\ &\cong R\R_{\Z_{2}}^{1}(\bar{\g}) / i^{*}R\R_{\Z_{2}}^{0}(\bar{\g})
\cong SR^{1}_{R}(\bar{\g}),
\end{split}\end{equation*}
using the commutative diagram above.
\end{proof}

In the real case, the representation group given by (\ref{eq:SR1}) is built from representations that are not only degree-shifted but also in fact $\bar{\g}$-modules rather than $\g$-modules.  We could not detect this subtlety when working with complex representations, since every complex $\g$-module is automatically a $\bar{\g}$-module and vice-versa (just multiply the action of the odd elements $X_{1}\in\g_{1}$ by $\sqrt{-1}$). We can now state the real version of Theorem \ref{exact-sequence}.

\begin{theorem}
	The degree-shifted real representation groups of a Lie superalgebra $\g$ form six-term periodic exact sequences \begin{equation}\label{eq:exact2}\begin{CD}
		R\R_{+}^{-n}(\g) @>\pi>> \SRR^{-n}(\g) @>\delta>> R\R_{\Z_{2}}^{-n}(\g) \\
		@Ai^{*}AA          @.         @VVfV \\
		R\R^{-n-1}_{\Z_{2}}(\g) @<\delta'<< \SRR^{n+1}(\bar{\g}) @<\pi'<< R\R_{+}^{-n-1}(\g)
	\end{CD}\end{equation}
\end{theorem}
\begin{proof}
This result follows immediately from Theorem \ref{exact-sequence} if we use the isomorphism $R\R_{Z_{2}}^{-n-1}(\g) \cong R\R_{0}^{-n}(\g)$ of equation (\ref{eq:real-grading-degree}), the isomorphism
$R\R_{+}^{-n-1}(\g) \cong R\R_{sc}^{-n}(\g)$ of Lemma \ref{restriction-equals-forgetful},
and Proposition \ref{restriction-forgetful-intertwine}.
\end{proof}

\begin{example}
For $\g = \R$ trivial, we have $\bar{\g} = \g$.  As before, the self-dual representation ideal is always
$\Rp^{-n}(\R) \cong \mathbb{Z}$ for all $n$, and the super representation
ring is the $KO$-theory of a point. We then have the exact
sequence
\begin{equation*}\begin{diagram}[tight,width=5em]
		\Z & \rTo^{\pi} & KO^{-n}(\mathrm{pt}) & \rTo^{\delta}& R\R^{-n}_{\Z_{2}}(\R)  \\
		\uTo>{i^{*}} & & & & \dTo>{f} \\
		R\R^{-n-1}_{\Z_{2}}(\R) & \lTo^{\delta} & KO^{n+1}(\mathrm{pt}) & \lTo^{\pi} &  \Z
\end{diagram}\end{equation*}
In fact, we obtain four distinct periodic exact sequences:
\begin{equation*}\begin{diagram}[tight,width=5em]
		\Z & \rTo^{0} & KO^{0}(\mathrm{pt}) \cong \Z & \rInto & R\R^{0}_{\Z_{2}}(\R) \cong \Z \oplus \Z \\
		\uTo>\cong & & & & \dOnto \\
		R\R^{-1}_{\Z_{2}}(\R) \cong \Z @<<< KO^{-7}(\mathrm{pt}) = 0 @<<< \Z
\end{diagram}\end{equation*}
\begin{equation*}\begin{diagram}[tight,width=5em]
		\Z & \rOnto & KO^{-1}(\mathrm{pt}) \cong \Z_{2} & \rTo^{0} & R\R^{-1}_{\Z_{2}}(\R) \cong \Z \\
		\uInto>{\times 2} & & & & \dTo>\cong \\
		R\R^{-2}_{\Z_{2}}(\R) \cong \Z @<<< KO^{-6}(\mathrm{pt}) = 0 @<<< \Z
\end{diagram}\end{equation*}
\begin{equation*}\begin{diagram}[tight,width=5em]
		\Z & \rOnto & KO^{-2}(\mathrm{pt}) \cong \Z_{2} & \rTo^{0} & R\R^{-2}_{\Z_{2}}(\R) \cong \Z \\
		\uInto>{\times 2} & & & & \dTo>\cong \\
		R\R^{-3}_{\Z_{2}}(\R) \cong \Z @<<< KO^{-5}(\mathrm{pt}) = 0 @<<< \Z
\end{diagram}\end{equation*}
\begin{equation*}\begin{diagram}[tight,width=5em]
		\Z & \rTo & KO^{-3}(\mathrm{pt}) = 0 & \rTo & R\R^{-3}_{\Z_{2}}(\R) \cong \Z \\
		\uOnto & & & & \dTo>\cong \\
		R\R^{-4}_{\Z_{2}}(\R) \cong \Z\oplus \Z & \lInto & KO^{-4}(\mathrm{pt}) \cong \Z & \lTo^{0}& \Z
\end{diagram}\end{equation*}
Due to the isomorphism $R\R_{\Z_{2}}^{-n}(\R) \cong R\R_{\Z_{2}}^{+n}(\R) \cong R\R_{\Z_{2}}^{n-8}(\R)$,
we see that the other four periodic exact sequences are obtained by
rotating each of these four by 180 degrees.
\end{example}


\end{document}